\documentclass{amsproc}
\usepackage{amssymb}

\newtheorem{predl}{Proposition}[section]
\newtheorem{defi}{Definition}[section]
\newtheorem{rem}{Remark}[section]

\newtheorem{lem}{Lemma}[section]
\newtheorem{theor}{Theorem}[section]
\newtheorem{example}{Example}[section]

\newcommand{\Si}{\Sigma}
\newcommand{\tr}{{\rm tr}}

\newcommand{\de}{\delta}
\newcommand{\al}{\alpha}

\newcommand{\be}{\beta}
\newcommand{\la}{\lambda}
\newcommand{\La}{\Lambda}

\newcommand{\ve}{\varepsilon}
\newcommand{\ep}{\epsilon}

\newcommand{\G}{\Gamma}
\newcommand{\ka}{\kappa}

\newcommand{\ga}{\gamma}

\newcommand{\om}{\omega}
\newcommand{\Om}{\Omega}

\def\clA{\mathcal{A}}
\def\clG{\mathcal{G}}
\def\clR{\mathcal{R}}
\def\clU{\mathcal{U}}

\def\clO{\mathcal{O}}
\def\clL{\mathcal{L}}
\def\clW{\mathcal{W}}
\def\clZ{\mathcal{Z}}

\def\we{\wedge}

\def\mC{{\mathbb C}}
\def\mZ{{\mathbb Z}}
\def\mR{{\mathbb R}}

\def\bfp{{\bf p}}

\def\bfM{{\bf M}}
\def\bfR{{\bf R}}

\newcommand{\thmat}[9]{\left(
\begin{array}{ccc}{#1}&{#2}&{#3}\\{#4}&{#5}&{#6}\\
{#7}&{#8}&{#9}
\end{array}\right)}
\newcommand{\beq}[1]{\begin{equation}\label{#1}}
\def\eq{\end{equation}}
\newcommand{\beqn}[1]{\begin{eqnarray}\label{#1}}
\newcommand{\eqn}{\end{eqnarray}}
\newcommand{\p}{\partial}
\newcommand{\di}{{\rm diag}}
\newcommand{\ti}{\tilde}
\newcommand{\oh}{\frac{1}{2}}

\newcommand{\GLN}{{\rm GL}(N,{\mC})}
\newcommand{\SLN}{{\rm SL}(N,{\mC})}
\newcommand{\SLt}{{\rm SL}(3,{\mC})}
\def\sln{{\rm sl}(N,\mC)}

\def\SL2{{\rm SL}(2,\mC)}

\newcommand{\ran}{\rangle}
\newcommand{\lan}{\langle}
\def\f1#1{\frac{1}{#1}}
\def\lb{\lfloor}
\def\rb{\rfloor}
\newcommand{\rar}{\rightarrow}

\newcommand{\bp}{\bar{\partial}}
\newcommand{\bz}{\bar{z}}
\newcommand{\bA}{\bar{A}}

\begin{document}

\vspace{0.3in}
\begin{flushright}
ITEP-TH-07/80\\
ESI- 1989
\end{flushright}
\vspace{0.3in}
\title{Lie Algebroids and generalized projective structures on Riemann surfaces}

\author{A.Levin}
\address{Institute of Oceanology, Moscow, Russia}
\email{alevin@wave.sio.rssi.ru}

\author{M.Olshanetsky}
\address{Institute of Theoretical and Experimental Physics, Moscow, Russia,}
\email{olshanet@itep.ru}
\thanks{We are grateful to A.Rosly and J.Grabowski for illuminating remarks.
The work were supported in part by  RFBR-06-02-17382, RFBR-06-01-92054-CE$_a$,\\
NSh-8065-2006.2.}



\maketitle
\begin{abstract}
The space of generalized projective structures on a Riemann surface $\Sigma$ of genus g with n marked points
is the affine space over the cotangent bundle to the space of SL(N)-opers.
It is a phase space of  $W_N$-gravity on $\Sigma\times\mathbb{R}$.
 This space is a generalization
of the space of projective structures on the Riemann surface.
We define the moduli space of  $W_N$-gravity
as a symplectic quotient with respect to the canonical action
of a special class of Lie algebroids.
They describe in particular the moduli space of deformations of complex structures
on the Riemann surface by differential operators of finite order, or equivalently,
by a quotient space of Volterra operators. We call these algebroids the Adler-Gelfand-Dikii (AGD)
algebroids, because they  are constructed by means of AGD bivector on the space of opers restricted on a circle.
The AGD-algebroids are particular case of Lie algebroids related to a Poisson sigma-model.
The moduli space of the generalized projective structure  can be described by cohomology of a BRST-complex.
\end{abstract}

\tableofcontents
\section {Introduction}
\setcounter{equation}{0}

The goal of this paper is to define a moduli space of generalized projective structures
on a Riemann surface $\Si_{g,n}$ of genus $g$ with $n$ marked points.
The standard projective structure is a pair of projective connection $T$ and the Beltrami
differential $\mu$. The projective connection is defined locally by the second order
differential operator $\p^2_z+T(z,\bz)$. It behaves as $(2,0)$-differential.
The Beltrami differential defines a deformation of complex structure on $\Si_{g,n}$ as
$\p_{\bz}$ $\to\,\p_{\bz}+\mu\p_z$. It is a $(-1,1)$-differential. The pair $(T,\mu)$
can be considered as coordinates in the affine space over the cotangent bundle to the space of projective connection. It is
an infinite dimensional symplectic space with the canonical form $\int_{\Si_{g,n}}DT\wedge D\mu$.
It has a field-theoretical interpretation as a phase space of $W_N$-gravity,
describing a topological field theory on $\Si_{g,n}\times\mR$ \cite{P,BFK,GLM}.
The chiral vector fields generate canonical transformations of this space.
The symplectic
quotient with respect to this action is a finite-dimensional space - the moduli space of
projective structure on   $\Si_{g,n}$.
The moduli space of complex structures is a part of
this space.

The projective connections have a higher order generalizations. It is a space of
$\SLN$-opers on   $\Si_{g,n}$ \cite{Tel,BD}. The dual to them variables define
generalized deformations of complex structures by $(-j,1)$-differentials $(j>1)$.
It can be equivalently
described as a quotient space of Volterra operators on  $\Si_{g,n}$.
These differentials describe
highest order integrals of motion in the Hitchin integrable systems \cite{Hi2}.
It is a base of the Hitchin fibration of the Higgs bundles and thereby their moduli parameterize the base.
In contrast with the space of projective structures ($W_2$-gravity)
the canonical transformations of this space generate Lie algebroids \cite{Ma,We}
rather than the Lie algebra of vector fields. The symplectic quotient with respect to this
action is a moduli space of generalized projective structures (the moduli space of $W_N$-gravity).
In particular, it describes the moduli space of generalized deformations of complex structures
by means of higher order differentials.

This situation can be generalized in the following way.
Remind that
 a Lie algebroid $\clA$ is a vector bundle over a space $M$ with Lie brackets defined
 on its sections $\G(\clA)$ and a bundle map (the anchor) to the vector fields
on the base $\de\,:\,\G(\clA)\to TM$. Let  ${\mathcal R}$ be an affine space  over
the cotangent bundle $T^*M$ (the principle homogeneous space).
It is a symplectic space with the canonical symplectic form.
Define a representation of  $\clA$ in the space of sections of $\clR$ in a such way
that together with the
anchor action they generate  canonical transformations of $\clR$.
 These transformations are classified by the first cohomology group
$H^1({\mathcal   A})$ of the algebroid.
 We call Lie algebroid equipped with the canonical representation in the
 space of sections
 ${\mathcal   R}$ {\sl the Hamiltonian algebroid} ${\mathcal   A}^H$.
The Hamiltonian algebroids are analogs of the Lie algebra of symplectic
vector fields.

The symplectic quotient of ${\mathcal   R}$ with respect to the canonical transformations
can be described by cohomology of the BRST operator of ${\mathcal   A}^H$
 \cite{BFV,HT}.
We prove that the BRST operator in this case has the same structure as
for the transformations by Lie algebras.

The general example of this construction is based on the Poisson sigma-model \cite{I,SS}.
We consider a one-dimensional Poisson sigma-model $X(t)\,:\,S^1\to M$,
where $M$ is a Poisson manifold.
The Lie algebroid in this case is a vector bundle over $\bfM=\{X(t)\}$.
The base $\bfM$ is Poisson with the Poisson brackets related to the Poisson
brackets on $M$. The space of sections of the Lie algebroid is $X^*(T^*M)$.
 The Lie brackets on  $X^*(T^*M)$  are defined by the Poisson brackets
on $\bfM$ \cite{GeD,Fu,MaMo}.
Though the Poisson structure exists only on
maps $S^1\to M$, in some important cases  the Lie algebroid structure (the Lie brackets
 and the anchor) can be continued  on the space of maps from a Riemann surface
 $\Sigma\supset S^1$ to $M$ and, moreover,
to the affine space $\bfR$ over $X^*(T^*M)$. The canonical transformations of the space of
smooth maps  $\bfR$ define a Hamiltonian algebroid.
It turns out that the first class constraints, generating these transformations are consistency
conditions for a linear system. Two of three linear equations define a deformation of operator
$\bp$ on $\Si$. This deformation depends on the Poisson bivector on $M$.
The symplectic quotient with respect to the canonical actions can be described by
the BRST construction. In particular, it defines
the moduli space of these deformations of the complex structure.

We apply this scheme to the space of opers on a Riemann surface $\Si_{g,n}$.
We start with the space $M_N(D)$ of $\GLN$-opers on a disk
$D\subset\Si$.
The space $M_N(D)$ being restricted on the boundary $S^1=\p D$ is a Poisson
space with respect to the Adler-Gelfand-Dikii (AGD) brackets \cite{Ad,GD}.
It allows us to define  a Lie algebroid (the AGD-algebroid) ${\mathcal   A}_N$ over $M_N(D)$.
 The Lie brackets on the space of sections
$\G({\mathcal   A}_N)$ and the anchor are derived from the AGD bivector.
The case $N=2$ corresponds to the projective structure on $D$
and the sections $\G({\mathcal   A}_N)$ is the Lie algebra of vector fields
with coefficients depending on the projective connection.
 The case $N>2$ is more involved and we deal with a genuine
Lie algebroids since differential operators with the principle symbol
of order two or more do not form a Lie algebra with respect to the standard commutator.
The AGD brackets define a new commutator on $\G({\mathcal   A}_N)$ that depends on
the projective connection and  higher spin fields.
According with the AGD construction the differential operators in $\G({\mathcal   A}_N)$
can be replaced by a quotient space of the Volterra operators on $D$.
This construction can be continued from $D$ to $\Si_{g,n}$.
We preserve the same notation for this algebroid and call it the AGD Lie algebroid
on $\Si_{g,n}$.
The space $M_N$ of opers on $\Si_{g,n}$ plays the role of the configuration space
of $W_N$-gravity \cite{P,BFK,GLM}. The whole phase space ${\bf   R}_N$
of $W_N$-gravity is an
affine space over the cotangent bundle to the space of opers $T^*M_N$.
  The canonical transformations of ${\bf R}_N$ are
sections of the Hamiltonian AGD-algebroid ${\mathcal   A}_N^H$ over $M_N$.
The symplectic quotient of the phase space
is the moduli space ${\mathcal   W}_N$ of the  $W_N$-gravity on $\Si_{g,n}$. Roughly
speaking, this space is
a combination of the moduli of $W_N$-deformations of complex structures
 and the spin $2\,,\ldots$,\,spin $N$ fields as the dual variables.
This moduli space can be described by the cohomology of
the BRST complex for the Hamiltonian algebroid.
 As it follows from
the general construction, the BRST operator has the same structure as in the Lie algebra case.
We consider in detail the simplest nontrivial case $N=3$. In this case it is possible
to describe explicitly the sections of the algebroid as the second order differential
operators, instead of Volterra operators.
It should be noted that the BRST operator for the $W_3$-algebras was constructed
in \cite{TM}. Here we construct the BRST operator for the different object
- the algebroid symmetries of $W_3$-gravity.  Another BRST description
of $W$-symmetries  was proposed in Ref.\cite{BL}.
We explain our formulae and the origin of the algebroid
by a special gauge procedure in the $\SLN$ Chern-Simons theory using an
approach developed in Ref.\cite{BFK}.

The paper is organized as follows. In next section we define
Lie algebroids, their representations, cohomolgy, the
Hamiltonian algebroids, and the BRST construction.
In Section 3 we use a Poisson sigma model to construct
Hamiltonian algebroids. In Section 4 we consider
two examples of our construction.
We analyze the moduli space of flat $\SLN$-bundles
and the moduli of projective structures on $\Si_{g,n}$.
A nontrivial example of this construction is $W_3$-gravity.
It is considered in detail in Section 5.
The general $W_N$ case  is analyzed in Section 6.

\section{Lie  algebroids and  groupoids}
\setcounter{equation}{0}

\subsection{Lie algebroids and groupoids}

We remind a definition brief description of Lie algebroids and Lie groupoids.
Details of this theory can be found in \cite{Ma,We}.

\begin{defi}
A {\em Lie algebroid} over a smooth manifold $M$ is a vector bundle
${\mathcal A}\rar M$ with
a Lie algebra structure on the space of its sections $\G({\mathcal A})$
 defined by the Lie brackets $\lb\ve_1,\ve_2\rb$,
$\ve_1,\ve_2\in\G({\mathcal A})$ and
a bundle map ({\em the anchor}) $\de :{\mathcal A}\to TM$, satisfying the following
conditions:
(i) For any $\ve_1,\ve_2\in\G({\mathcal A})$
\beq{5.1}
[\de_{\ve_1},\de_{\ve_2}]=\de_{\lb\ve_1,\ve_2\rb}\,,
\end{equation}
(ii) For any $\ve_1,\ve_2\in\G({\mathcal A})$ and $f\in C^\infty(M)$
\beq{5.2}
\lb\ve_1,f\ve_2\rb=f\lb\ve_1,\ve_2\rb + (\de_{\ve_1} f)\ve_2\,.
\end{equation}
\end{defi}
In other words, the anchor  defines a representation of $\G({\mathcal A})$
in the Lie algebra of vector fields
on $M$. The second condition is the Leibnitz rule with respect to the multiplication of
 sections by smooth functions.

Let $\{e^j(x)\}$ be a basis of sections. Then  the brackets
are defined by the structure functions $f^{jk}_i(x)$ of the algebroid
\beq{5.1b}
\lb e^j,e^k\rb=f^{jk}_i(x)e^i\,,~~x\in M\,.
\end{equation}
Using the Jacobi identity for the brackets $\lb\,,\,\rb$, we find
\beq{5.5}
f^{jk}_i(x)f_n^{im}(x)+\de_{e^m}f^{jk}_n(x)+{\rm c.p.}(j,k,m)=0\,.
\footnote{The sums over repeated indices are understood throughout the paper,
and ${\rm c.p.}(j,k,m)$ means the cycle permutation.}
\end{equation}
If the anchor is trivial, then ${\mathcal   A}$ is just a bundle of Lie algebras.\\

There exists the global object - the {\em  Lie groupoid} \cite{GF}.

\begin{defi}
A Lie groupoid $G$ over a manifold $M$
is a pair of smooth manifolds $(G,M)$,  two
smooth mappings $l,r~: ~G\to M$ and a partially defined smooth binary operation (the
product)
$(g,h)\mapsto g\cdot h $ satisfying the following conditions:\\
(i) It is defined when $l(h)=r(g)$. \\
(ii) It is associative: $(g\cdot h)\cdot k=g\cdot (h\cdot k)$
whenever the products are defined.\\
(iii) For any $g\in G$ there exist the left and right identity elements $l_g$ and
$r_g$ in $G$ such that $l_g \cdot g=g\cdot r_g=g$.\\
(iv) Each $g$ has an inverse $g^{-1}$ such that $g\cdot g^{-1}=l_g$ and
$g^{-1}\cdot g=r_g$.\\
\end{defi}
We denote an element of $g\in G$ by the triple $\lan\lan x|g|y\ran\ran$, where
$x=l(g),~y=r(g)$.
Then  the product $g\cdot h $ is
$$
g\cdot h= \lan\lan x|g\cdot h|z\ran\ran=\lan\lan x|g|y\ran\ran\lan\lan y|h|z\ran\ran\,.
$$
An orbit of the groupoid in the base $M$ is defined as an equivalence
$x\sim y$ if $x=l(g),~y=r(g)$ for some $g\in G$. The isotropy subgroup $G_x$ for $x\in M$ is
defined as
$$
G_x=\{g\in G~|~l(g)=x=r(g)\}=\{\lan\lan x|g|x\ran\ran\}\,.
$$

The Lie algebroid is a infinitesimal  version of the Lie groupoid. The anchor is determined
in
terms of the multiplication law.
The problem of integration of a Lie algebroid to a Lie groupoid is treated
in Ref.\,\cite{GF}.

\subsection{Representations and cohomology of Lie algebroids}

The definition of algebroids representations is rather evident:
\begin{defi}
A vector bundle representation (VBR) $(\rho, {\mathcal   M})$ of a Lie algebroid
$\mathcal   A$ over $M$
is a vector bundle $\mathcal   M$ over $M$ and a map
 $\rho$ from $ {\mathcal   A}$ to the bundle of differential operators in $\mathcal   M$
  $\mbox{\it Diff}^{~\le 1}({\mathcal   M}, {\mathcal   M})$ of the order
 less or equal to $1$,
 such that:\\
(i) the principle symbol of $\rho(\ve)$ is a scalar equal to the anchor of $\ve$:
\beq{symb}
{\rm Symb}(\rho(\ve))={\rm Id}_{\mathcal   M}\de_\ve\,,
\end{equation}
(ii) for any $\ve_1,\ve_2\in\G({\mathcal   A})$
\beq{rep}
[\rho({\ve_1}),\rho({\ve_2})]=\rho({\lb\ve_1\ve_2\rb})\,,
\end{equation}
 where  the l.h.s. denotes
the commutator of differential operators.
\end{defi}

For example, the trivial bundle is a VBR representation (the map $\rho$
is the anchor map $\delta$).

Consider a small ball $U_{\alpha} \subset M$ with local coordinates
$x=(x^1,\ldots,x^n)$. Then the anchor can be written as
\beq{5.15i}
\de_{e^j}=b^{ja}(x)\frac{\p}{\p x^a}=\lan b^j|\p\ran\,,~~(a=1,\ldots,\dim M)   \,.
\footnote{The brackets $\lan~|~\ran$ mean summations over repeating indices.}
\end{equation}
Let $\xi=\xi=\xi^a\p_a\,,$ $(\p_a=\frac{\p}{\p x^a})$
 be a section of the tangent bundle $TM$. Then the VBR on
 $TM$ takes
the form
\beq{5.15j}
{\rho}(e^j)\xi=\lan b^j|\p \xi\ran   -\lan\p b^j|\xi\ran=[\de_{e^j},\xi]\,.
\end{equation}

Similarly, the VBR on a section $p=p_adx^a=\lan p|dx\ran$ of $T^*M$ is
\beq{5.16i}
\rho(e^j)p=\lan b^j|\p\ran\lan p|dx\ran+\lan dx|\lan\p(b^j)p\ran\ran
={\mathcal L}_{\de_{e^j}}p \,,
\end{equation}
where ${\mathcal L}$ is the Lie derivative along the vector field $\de_e^j$.
We omit a more general definition of the sheaf representation.

Now we define  cohomology groups of  algebroids.
First, we consider the case of contractible base $M$.
Let ${\mathcal   A}^*$ be a bundle over $M$ dual to ${\mathcal   A}$.
 Consider the bundle of graded commutative algebras
$\wedge^\bullet{\mathcal   A}^*$.
The space $\G(M, \wedge^\bullet {\mathcal   A}^*)$
is generated by the sections $\eta_k$:
$\lan\eta_j|e^k\ran   =\de_j^k$. It is a graded algebra
$$
\G(M, \wedge^\bullet  {\mathcal   A}^*)=\oplus {\mathcal   A}^*_n\,,~~
{\mathcal   A}^*_n=\{c_n(x)=
\f1{n!}c_{j_1,\ldots,j_n}(x)\eta_{j_1}\ldots\eta_{j_n}\,,~x\in U\}\,.
$$
$$
c_{j_1,\ldots,j_n}(x)\in\mathcal{O}(M)\,.
$$

 Define the Cartan-Eilenberg operators ``dual'' to the  brackets
 $\lb ,\rb$
$$
sc_n(x;e^{1},\ldots,e^{n},e^{n+1})=\sum_i
(-1)^{i-1}\de_{e^i}c_n(x;e^{1},\ldots,\hat{e}^{i},\ldots e^{n})-
$$
\beq{5.6}
-\sum_{j<i}
(-1)^{i+j}c_n(x;\lb e^{i},e^j\rb,\ldots,\hat{e}^{j},\ldots,\hat{e}^{i},\ldots,e^{n})\,,
\end{equation}
where
$$
\de_{e^i}c_n(x)=\lan \p c_n(x)|b^i(x)\ran\,.
$$
It follows from (\ref{5.1}) and  (\ref{5.5}) that $s^2=0$.
Thus, $s$ defines a complex of bundles ${\mathcal   A}^*\to \wedge^2{\mathcal   A}^*\to
\cdots$.

The cohomology groups $H^k(\mathcal{A},\mathcal{O}(M))$  of this complex are called {\em the cohomology groups of
algebroid with trivial coefficients}.
This complex is  a part of the BRST complex described below.

 The action of the coboundary operator $s$ takes the following form
 on the low cochains:
\beq{5.7}
sc(x;\ve)=\de_\ve c(x)\,,
\end{equation}
\beq{5.8}
sc(x;\ve_1,\ve_2)=\de_{\ve_1}c(x;\ve_2)-\de_{\ve_2}c(x;\ve_1)-
c(x;\lb \ve_1,\ve_2\rb)\,,
\end{equation}
\beq{5.9}
sc(x;\ve_1,\ve_2,\ve_3)=\de_{\ve_1}c(x;\ve_2,\ve_3)-
\de_{\ve_2}c(x;\ve_1,\ve_3)
\end{equation}
$$
+\de_{\ve_3}c(x;\ve_1,\ve_2)
-c(x;\lb \ve_1,\ve_2\rb,\ve_3)
+c(x;\lb \ve_1,\ve_3\rb,\ve_2)
-c(x;\lb \ve_2,\ve_3\rb,\ve_1)\,.
$$
It follows from (\ref{5.7}) that $H^0({\mathcal   A},\clO(M))$  is isomorphic to the
invariants in the space $\mathcal{O} (M)$.

The next cohomology group
$H^1({\mathcal   A},\clO(M))$ is responsible for the shift of the anchor action on
$\clO(M)$
$\de_\ve f(x)=\lan\de_\ve x|\p f(x)\ran)$
\beq{5.11}
\hat{\de}_\ve f(x)= \de_\ve f(x)+c(x;\ve)\,, ~~sc(x;\ve)=0\,.
\end{equation}
Then due to (\ref{5.8}), this action is
 consistent with the defining anchor property (\ref{5.1}).

We modify the action (\ref{5.11}). Let $\Psi=\exp f\in\mathcal{O}^* (M) $.
Then the action on $\Psi$ takes the form
\beq{5.11a}
\ti{\de}_{\ve}\Psi(x)=\left(\de_\ve+c(x;\ve)\right)\Psi(x)
\end{equation}
satisfies the the anchor property (\ref{5.1})
$[\ti{\de}_{\ve_1}\ti{\de}_{\ve_2}]=\ti{\de}_{[\ve_1,\ve_2]}$.

If $M$ is not contractible the definition of cohomology group is more
complicated. We sketch the \^{C}ech version of it.
Choose an acyclic covering $U_\alpha$. Consider the \^{C}ech
complex with coefficients in $\bigwedge{}^{\bullet}({\mathcal   A}^{*})$
corresponding to this covering:
$$
\bigoplus \Gamma (U_\alpha ,\bigwedge{}^{\bullet}({\mathcal A}^{*}))
\stackrel{d^C}{\longrightarrow}
\bigoplus\Gamma (U_{\alpha\beta },\bigwedge{}^{\bullet}({\mathcal A}^{*}))
\stackrel{d^C}{\longrightarrow}\cdots
$$
The \^{C}ech differential $d^C$ commutes with  the Cartan-Eilenberg operator
$s$, and cohomology of algebroid are cohomology of normalization of this
bicomplex~:
$$\bigoplus \Gamma (U_\alpha ,{\mathcal A}^{*}_0)
\stackrel{d^C,s}{\longrightarrow}
\bigoplus\Gamma (U_{\alpha\beta },{\mathcal A}^{*}_0)
\oplus
\bigoplus \Gamma (U_\alpha ,{\mathcal A}^{*}_1)
\stackrel{
\left(\begin{array}{ccc}
d^C&s&\\
&-d^C&s
\end{array}\right)
}{\longrightarrow}$$
$$
\bigoplus\Gamma (U_{\alpha\beta\gamma },{\mathcal A}^{*}_0)
\oplus
\bigoplus \Gamma (U_{\alpha\beta} ,{\mathcal A}^{*}_1)
\oplus
\bigoplus \Gamma (U_\alpha ,{\mathcal A}^{*}_2)\longrightarrow \cdots\,.
$$
The cochains
 $c^{i,j}\in\bigoplus_{\alpha_{1}\alpha_{2}\cdots \alpha_j}
\Gamma (U_{\alpha_{1}\alpha_{2}\cdots \alpha_j} ,{\mathcal A}^{*}_i)$ are bigraded.
The differential maps $c^{i,j}$ to $(-1)^jd^C\,c^{i,j}+sc^{i,j}$,
 has type $(i,j+1)$ for $(-1)^jd^C\,c^{i,j}$ and $(i+1,j)$ for $sc^{i,j}$.

Again, the group $H^0({\mathcal A},\clO(M))$  is isomorphic to the
invariants in the whole space $\mathcal{O} (M)$.

Consider the next group $H^{(1)}({\mathcal A},\clO(M))$.
It has  two components\\ $(c_\alpha(x,\ve),c_{\alpha\beta}(x))$.
They are characterized by the following conditions (see (\ref{5.8}))
 $$
c_\alpha(x;\lb \ve_1,\ve_2\rb)=
\de_{\ve_1}c_\alpha(x;\ve_2)-\de_{\ve_2}c_\alpha(x;\ve_1)\,,
$$
\beq{ch1}
  \de_{\ve}c_{\alpha\beta}(x)=-c_\alpha(x;\ve)+c_\beta(x;\ve)\,,
\end{equation}
\beq{ch2}
 c_{\alpha\gamma}(x)=c_{\alpha\beta}(x)+c_{\beta\gamma}(x)\,.
\end{equation}

The group  $H^{(2)}({\mathcal   A},\clO(M))$ is responsible for the central
extension of the the brackets on $\G({\mathcal   A})$. Let  $c(x;\ve_1,\ve_2)$  be a
 two-cocycles. Then
\beq{5.9a}
\lb(\ve_1,0),(\ve_2,0)\rb_{c.e.}=
(\lb\ve_1,\ve_2\rb,c (x;\ve_1,\ve_2))\,.
\end{equation}
The cocycle condition (\ref{5.9}) means that the new brackets
$\lb~,~\rb_{c.e.}$ satisfies  (\ref{5.5}).
The exact cocycles lead to the split extensions.
There is an obstacle to this continuations in $H^{(3)}({\mathcal   A},\clO(M))$.
We do not dwell on this point.

\begin{example}
Consider a flag variety $Fl_N=G/B$, where $G=\SLN$ and $B$ is the lower Borel subgroup.
The flag variety is a base of a $\sln$ algebroid.
 In particular,  sl$(2,\mC)$ anchor acts on $Fl_2\sim \mC P^1$
in a neighborhood $\clU_+$ of $z=0$ by the vector fields
$$
\de_e=\p_z\,,~~\de_h= -2z\p_z\,,~~\de_f=-z^2\p_z\,.
$$
The one-cocycle representing $H^1({\rm sl}(2,\mC),\clO(\mC P^1))$
$$
c(z,e)=0\,,~~c(z,h)=\nu\,,~~c(z,f)=\nu z
$$
defines the extension of the anchor action
$\hat{\de}_\ve$ (\ref{5.11}).
Let $w$ be a local coordinate in a neighborhood $\clU_-$ of $z=\infty$.
Then
$$
\de_e=-w^2\p_w\,,~~\de_h= 2w\p_w\,,~~\de_f=\p_w\,.
$$
and
$$
c(w,e)=\nu w\,,~~c(w,h)=-\nu\,,~~c(w,f)=0\,.
$$
The another component of  $H^1({\rm sl}(2,\mC),\clO(\mC P^1))$
satisfying (\ref{ch1}) is the cocycle
\beq{coc}
c_{+-}(z)=\nu\log z\,.
\end{equation}
\end{example}

 \begin{example}
Consider the $\mC^*$ bundle over $\mC^2\setminus 0=\{(z_1,z_2)\}$.
Define  the anchor map
\beq{c2}
\de_\ve=\ve(z_1\p_{z_1}+z_2\p_{z_2})\,.
\end{equation}
 According with (\ref{5.11}) it can be extended as
\beq{cp1}
\hat{\de}_\ve=\ve(z_1\p_{z_1}+z_2\p_{z_2}-\nu)\,,
\end{equation}
because $\ve\nu\in\mC$ represents a non-trivial one-cocycle.
\end{example}

\bigskip





\subsection{Affine spaces over cotangent bundles}

We shall consider Hamiltonian algebroids over cotangent bundles.
The cotangent bundles
are a special class of symplectic manifolds. There exist a generalization of
 cotangent bundles, that we include in our scheme.
 It is  affine spaces over a cotangent bundles we are going to
define.

Let $V$ be a vector space and $\clR$ is a manifold with an action of
$V$ on $\clR$
$$
\clR\times V\rar \clR:(x,v)\to x+v\in \clR\,.
$$
\begin{defi}
The manifold ${\mathcal   R}$ is an {\sl affine space over $V$} (a {\em  principle
homogeneous
space over} $V$) if the action of $V$
on  ${\mathcal   R}$ is transitive and free. We denote it ${\mathcal   R}/ V$.
\end{defi}
In other words, for any pair $x_1,x_2\in{\mathcal   R}$ there exists $v\in  V$ such that
$x_1+v=x_2$, and $x+v\neq x$ if $v\neq 0$.

This construction is generalized to vector bundles. Let $E$ be a vector bundle
 over $ M$ and $\G(E)$ is the linear space of its sections.
\begin{defi}
 An affine space ${\mathcal   R}/E$ over  $E$
 is a bundle over $M$ with the space of
 sections $\G({\mathcal   R})$ defined as the affine space over $\G(E)$.
\end{defi}

Consider a cotangent bundle  $T^*M$ and the corresponding
affine space ${\mathcal   R}/T^*M$. Let  $\xi_\al$ be a section
of ${\mathcal   R}/T^*M$ over a contractible set  $\mathcal{U}_{\al}\subset M$.
It can be identified with a section of $T^*\mathcal{U}_{\al}$.
 To define ${\mathcal R}$ over $M$ consider
an intersection $\mathcal{U}_{\al\be}=\mathcal{U}_\al\cap\mathcal{U}_\be$.
We assume that
local sections are related as
\beq{ral}
\xi_\al=\xi_\be+\varsigma_{\al\be}\,,
\end{equation}
where $\varsigma_{\al\be}\in\G(\clZ^{(1)}(\mathcal{U}_{\al\be}))$
and
$$
\clZ^{(1)}(\mathcal{U}_{\al\be})=\{\varsigma_{\al\be}\in\Om^{(1)}(\mathcal{U}_{\al\be})\,
|\,d\varsigma_{\al\be}=0\}\,.
$$
We have $\varsigma_{\al\be}=-\varsigma_{\be\al}$ and  on the intersection
$\mathcal{U}_{\al\be\ga}=\mathcal{U}_\al\cap\mathcal{U}_\be\cap\mathcal{U}_\ga$
$$
\varsigma_{\al\be}+\varsigma_{\be\ga}+\varsigma_{\ga\al}=0\,.
$$
If $\varsigma_{\al\be}=p_\al-p_\be$ then $\clR$ can be continued from
$\clU_\al$ to $\clU_\be$ as $T^*\mathcal{U}_{\be}$.
Therefore, the non-equivalent affine spaces over $T^*M$ are classified by the
elements $H^1(M,\clZ^{(1)})$.

\bigskip


An affine space ${\mathcal   R}/T^*M$ is
a symplectic space with the canonical form $\om$.
Locally  on $\mathcal{U}_{\al}$, $\om=\lan  d\xi_\al\wedge dx_\al\ran$.
In fact, the symplectic form is well defined globally, because
 the transition forms $\varsigma_{\al\be}$ are closed.
But in contrast with $T^*M$,
$\om$ is not exact, since $\xi_\al dx_\al$ is defined only locally.
The symplectic form on $T^*M$ has vanishing cohomology class
$[\om_{T^*M}]=0$, while the cohomology class of $[\om_{\clR}]$ is an
nontrivial element in $H^2(M,\mC)$.

\bigskip
\texttt{Example} {\bf 2.3}

The cotangent bundle $T^*Fl_N$ to the flag variety \\
$Fl_N=\SLN/B$,
 without the null section,
can be identified with the coadjoint orbit passing through
the Jordanian matrix $\sum_{j=1}^{N-1}E_{j,j+1}$. Consider
in particular $Fl_2\sim\mC P^1$.
The symplectic form on the cell $\clU_+=\{|z|<\infty\}\subset \mC P^1$
 is $\om=dp_+\wedge dz_+$.
The cotangent bundle $T^*\clU_+$ can be represented by the matrix
$$
\left(
    \begin{array}{cc}
      -z_+p_+ & p_+\\
     -z^2_+p_+  & z_+p_+ \\
    \end{array}
  \right)\,.
$$

 On $\clU_-=\{|z|>0\}$ the form $\om$ is
 $dp_-\wedge dz_-\,$ $\,(z_-=\f1{z_+})$. On the intersection the transition
 form vanishes $p_+dz_+=p_-dz_-$.

 The affine space ${\mathcal R}/T^*\mC P^1$
is a generic coadjoint orbit $\mathcal{O}_\nu $
passing through $\di(\nu,-\nu)$. Over  $\clU_+$ the coadjoint orbit
has the parametrization $(\xi_+,z_+)$
$$
\left(
    \begin{array}{cc}
      -z_+\xi_++\oh\nu & \xi_+\\
     -z^2_+\xi_++\nu z_+  & z_+\xi_+ -\oh\nu\\
    \end{array}
  \right)
$$
with the form $d\xi_+\we dz_+$. On the cell $\clU_-$ the form is
$d\xi_-\we dz_-$. The transition form is represented by the non-trivial
cocycle from $H^1(\mC P^1,T^*Fl_2)$
\beq{tra}
\xi_-dz_--\xi_+dz_+=\nu\, d(\log z_-)
\end{equation}
(compare with (\ref{coc})).

\texttt{Example} {\bf 2.4}

The basic example, though  infinite-dimensional, is the affine space over
the antiHiggs bundles.
\footnote{We use the antiHiggs bundles instead of the standard Higgs bundles for
reasons, that will become clear in Sect. 4.}
The antiHiggs bundle ${\mathcal   H}_N(\Si)$ is a cotangent bundle to the space
of holomorphic connections $M=\{\nabla^{(1,0)}=\p+A \}$
in a trivial vector bundle of rank $N$ over a complex curve
$\Si$. The cotangent vector (the antiHiggs field) is $\sln$ valued $(0,1)$-form
$\bar{\Phi}$. The symplectic form on ${\mathcal   H}_N(\Si)$ is
$-\int_\Si\tr(D\bar{\Phi}\wedge DA)$.
An example of the affine space ${\mathcal   R}/{\mathcal   H}_N(\Si)$
is the space of the $\ka$-connections $\{(\kappa\bp+\bA,\p+A)\}\,$ $(\ka\in\mC)$ with the
symplectic form $\int_\Si\tr(DA\wedge D\bA)$.


\subsection{Lie algebroids with representations in affine spaces}

 Let  ${\mathcal   A}$ be a Lie algebroid  over $M$
with the anchor $\de$. We shall define a representation  $\rho$  of  ${\mathcal   A}$
acting on sections of the affine space ${\mathcal   R}/T^*M$ in such a way that
 prove that $\rho$ combined with  $\de$ are hamiltonian vector fields on ${\mathcal   R}$
 with respect to $\om=\lan  d\xi|dx\ran$.
  It means that
 $$
 i_\ve\om=dh_\ve\,,~~\de_\ve x=\{h_\ve,x\}\,,~~\rho_\ve \xi=\{h_\ve,\xi\}\,.
 $$

\begin{lem}
The anchor action (\ref{5.15i}) of a Lie algebroid ${\mathcal   A}$ over $M$
can be lifted to a hamiltonian action on ${\mathcal   R}/T^*M$
with
 \beq{5.17}
 h_\ve=\lan \de_\ve |\xi\ran   +c(x;\ve)\,,
\end{equation}
where  $c(x;\ve)\in H^1(\clA,\clO(M))$.
\end{lem}


{\sl Proof}.\\
Define first the VBR of $\clA$ on sections of $\clR$.
Consider a contractible set $\mathcal{U}_\al\subset M$
with local coordinates $x_\al=(x^1_\al,\ldots,x^d_\al)$
 The anchor  has the form
\beq{5.15a}
\de_{\ve}x^a_\al= V_\al^a\,,~~  (\,\de_{e^j}=\sum_a V^{j,a}\p_a=\lan V^j|\p\ran\,)
\end{equation}
Let $\xi_\al\in\G(\mathcal{R}_\al/T^*\mathcal{U}_\al)$.
Define  an action $\rho$  on the sections
\beq{5.16a}
\rho_{\ve}\xi_\al=
\clL_{\ve}\xi_\al+dc_\al(x;\ve)\,,
\end{equation}
where $\clL_{\ve}$ is the Lie derivative and $c_\al(x;\ve)$ represents an element from $H^1(\clA,\clO(M))$.
 Since $c_\al(x;\ve)$ is a cocycle,
the action (\ref{5.16a}) satisfies (\ref{rep}).
It also satisfies (\ref{symb}). Therefore, the action (\ref{5.16a}) is a VBR.
The last term in (\ref{5.16a}) is responsible for the passage from $T^*\clU_\al$ to
the affine space ${\mathcal   R}_\al$, otherwise $\xi_\al$
 is transformed as a cotangent vector
(see (\ref{5.16i})).

To define the VBR globally we prove that
 on the intersection
$\clU_{\al\be}=\clU_\al\cap\clU_\be$ we have
\beq{ide}
\rho_{\ve}(\xi_\al-\xi_\be)=\clL_\ve(\varsigma_{\al\be}+dc_{\al\be})\,,
\end{equation}
where $\varsigma_{\al\be}$ is a closed one-form representing an element from
$H^1(M,\clZ^{(1)})\,$ and
$c_{\al\be}\in H^1(M,\clO(M))$ (see (\ref{ch2})).
To prove (\ref{ide})
we apply  (\ref{5.16a}) to its left hand side
$$
\rho_{\ve}(\xi_\al-\xi_\be)=\clL_\ve(\xi_\al-\xi_\be)+d(c_\al(x;\ve)-
c_\be(x;\ve))\,.
$$
Since $d(\xi_\al-\xi_\be)=d\varsigma_{\al\be}=0\,$, $\clL_\ve(\xi_\al-\xi_\be)=d\imath_\ve\varsigma_{\al\be}$.
Then using (\ref{ch1}) we come to (\ref{ide}).
The action (\ref{ide}) means that
$$
\rho_{\ve}\varsigma_{\al\be}=\clL_\ve(\varsigma_{\al\be}+dc_{\al\be})
$$
is the Lie derivative of a closed one-form.
It allows us to define the VBR on sections of ${\mathcal
R}/T^*M$.

The direct calculations show that $\de_\ve$ (\ref{5.15a}) and $\rho_\ve$
(\ref{5.16a}) are hamiltonian vector fields $\{h_\ve\,,~\}$ with
 $h_\ve$  (\ref{5.17}). The  Hamiltonians have the linear dependence on "momenta".
The corresponding  Hamiltonians (\ref{5.17}) have the linear dependence on "momenta".
The exact cocycle $c(x;\ve)=\de_{\ve}f(x)$
shifts $\xi\,$ in the Hamiltonian
$h_\ve=\lan \xi+d f(x)|\de_\ve)\ran$.
Thus, all nonequivalent
lifts of anchors from $M$ to ${\mathcal   R}/T^*M$ are in one-to-one correspondence
with $H^1({\mathcal   A},\clO(M))$.
$\Box$

\begin{rem}
There exists the map (see (\ref{5.17}))
\beq{5.17aa}
\G(\clA)\,\to\,\clO(\clR)\,,~~~~(\ve\,\to\,h_\ve=\lan \de_\ve |\xi\ran+c(x;\ve))\,.
\eq
Due to (\ref{5.1}) and the cocycle property of $c(\ve,x)$ it is the Lie algebras map
\beq{5.12}
\{h_{\ve_1},h_{\ve_2}\}=h_{[\ve_1,\ve_2]}\,.
\end{equation}
The map to the hamiltonian vector fields is the bundle map
$$
f\ve\,\to f\{h_\ve,~\}\,,
$$
because the ring of functions is defined on the base $M$ and thereby is Poisson commutative.
\end{rem}
This remark suggests the following definition.
\begin{defi}
We call
the Lie algebroid $\clA$ over $M$ equipped with the VBR in the sections of
$\clR/T^*M$ the {\it Hamiltonian algebroid} $\clA^H$. The anchor of $\clA^H$
is the map (\ref{5.17aa}).
\end{defi}

Then the Jacobi identity for the Hamiltonian algebroids assumes the form
\beq{5.14}
f^{jk}_i(x)f_n^{im}(x)+\{h_{e^m},f^{jk}_n(x)\}+c.p.(j,k,m)=0\,.
\end{equation}

\bigskip

\texttt{Example} {\bf 2.5}

Consider the  action of $\mC^*\,:$ on
$\mC^2\setminus \{0\}\,$, $\,(z_a\to \la z_a)\,,$
$\,(\la\neq 0)\,.$
Infinitesimally, we have   (\ref{c2}) (see Example 2.2).
  The quotient $(\mC^2\setminus \{0\})/\mC^*$
is isomorphic to $\mC P^1\sim Fl_2$.
We lift the infinitesimal action (\ref{c2})
to the  cotangent bundle $T^*(\mC^2\setminus \{0\})$.
The cotangent bundle $T^*(\mC^2\setminus \{0\})$ is equipped with
the canonical symplectic form
$$
\om=dp_1\wedge dz_1+dp_2\wedge dz_2\,.
$$
It is invariant under the action of $\mC^*\,:$ $\,z_a\to(\exp\ve) z_a\,,$
$p_a\to (\exp-\ve) p_a$.
The generating Hamiltonian $h_\ve$ of this transformation $(\imath_\ve \om=d\mu^*)$
has the form
$$
h_\ve=-\ve(p_1z_1+p_2z_2)\,.
$$
It defines the Hamiltonian algebroid $\clA^H$.

 Consider the symplectic quotient of  $T^*(\mC^2\setminus 0)$ with respect to the $\mC^*$ action generating by the
shifted Hamiltonian by the cocycle $\ve\nu$
$$
h_\ve\to h_\ve-\ve\nu\,.
$$
It defines  the moment map constraint
$$
p_1z_1+p_2z_2=\nu\,.
$$
For $z_1\neq 0$ one can fix the gauge $z_1=1$ and find from the moment constraint $p_1=-p_2z_2+\nu$. In this case the symplectic quotient
is isomorphic to $T^*\mC$
with coordinates $(p_2,z_2)$ and the form $dp_2\we dz_2$.
Similarly, for $z_2\neq 0$ one can take $z_2=1\,,$  $\,p_2=-p_1z_1+\nu$
and the canonical coordinates $(p_1,z_1)$ on the symplectic quotient $T^*\mC$.
On the intersection we come to the relation
$$
p_1dz_1=p_2dz_2+d\log z_1^\nu\,.
$$
Comparing with (\ref{tra}) we conclude that for a nonzero value of the moment map
$\nu\neq 0 $ the symplectic quotient is a generic coadjoint orbit, otherwise for
$\nu=0$ we come to $T^*\mC P^1$.

\subsection{Reduced phase space and its BRST description}

 Let $e^j$ be a basis of sections in $\G({\mathcal   A})$.
Then the Hamiltonians (\ref{5.17}) can be represented in the form
$h^j=\lan e^j|F(x)\ran   $, where $F(x)\in\G({\mathcal   A}^{*})$ defines the {\em moment
map}
$$
F\,: {\mathcal   R}\rar\G({\mathcal   A}^{*})\,.
$$
The coadjoint action ${\rm ad}^*_{\ve}$ in $\G(({\mathcal   A}^{H}))$ is defined
in the standard way
$$
\lan \lb\ve, e^j\rb|F(x)\ran   =\lan e^j|{\rm ad}^*_{\ve}F(x)\ran \,.
$$
The zero-valued moment  $F(x)=0$ is preserved by the groupoid
coadjoint action  $G$ generated by ad$^*_\ve$.
The moment constraints $F(x)=0$ generate
canonical algebroid action on ${\mathcal R}$.
 The reduced phase space is defined as the quotient
$$
{\mathcal   R}^{red}={\mathcal   R}//G:=\{x\in {\mathcal   R}|(F(x)=0)/G\}\,,
$$
 In other words, ${\mathcal R}^{red}$
is the set of orbits of $G$ on the constraint surface\\
 $F(x)=0$.

The BRST approach allows us to go around the reduction procedure by introducing
additional fields (the ghosts).
We shall construct the BRST complex for ${\mathcal   A}^{H}$ in a similar way as the
Cartan-Eilenberg complex for the Lie algebroid ${\mathcal   A}$. The BRST complex
is endowed with a Poisson structure and in this way it has the form of the
quasi-classical limit of the Cartan-Eilenberg complex.

 The sections
  $\eta\in\G(({\mathcal   A})^*)$ are the anti-commuting (odd) fields called
 {\em the ghosts}.
We preserve the notation for the Hamiltonians in terms of the ghosts
$h=\lan \eta|F(x)\ran$, where $\{\eta_j\}$ is a basis in
$\G(({\mathcal   A})^*)$.
 Introduce another type of odd variables ({\em the ghost momenta})
${\mathcal   P}^j\in \G({\mathcal A}^{H})$ dual to the ghosts $\eta_k$.
 We attribute the ghost
 number one to the ghost fields gh$(\eta)=1$, minus one to the ghost momenta
gh$({\mathcal   P})=-1$ and gh$(x)=0$ for $x\in {\mathcal   R}$.
Define
the Poisson brackets in addition to the Poisson structure on ${\mathcal R}$
\beq{5.19}
\{\eta_j,{\mathcal   P}^k\}=\de_j^k\,,~~\{\eta_j,x\}=\{{\mathcal   P}^k,x\}=0\,.
\end{equation}
All fields are incorporated in the graded Poisson superalgebra
$$
{\mathcal BFV}=\left(
\G(\wedge^\bullet  (({\mathcal   A})^*\oplus{\mathcal   A})
\right)\otimes \clO({\mathcal   R})
=\G(\wedge^\bullet ({\mathcal   A})^*)
\otimes\G(\wedge^\bullet {\mathcal   A})\otimes \clO({\mathcal   R})\,.
$$
({\em the Batalin-Fradkin-Volkovitsky (${\mathcal BFV}$) algebra}).

There exists a nilpotent operator $Q$ on the ${\mathcal BFV}$ algebra $Q^2=0,~gh(Q)=1$
({\em the BRST operator}) transforming it into the BRST complex.
The cohomology of the BRST complex give rise to the structure of
the classical reduced phase space ${\mathcal R}^{red}$. Namely, $H^0(Q)$ is
the space of invariants with respect to the symplectic action of $G$ and
in this way can be identified with the classical observables.

Suppose that the action of $Q$ has the hamiltonian form:
$$
Q\psi=\{\psi,\Om\},~~\psi,\Om\in {\mathcal BFV}\,.
$$
Due to the Jacobi identity the nilpotency of $Q$
is equivalent to
$\{\Om,\Om\}=0$.
Since $\Om$ is odd, the brackets are symmetric.
 For a generic Hamiltonian system with the first class constraints $\Om$
 can be represented as the expansion \cite{HT}
$$
\Om=h_\eta+\oh\lan \lb\eta,\eta'\rb|{\mathcal   P}\ran   +...\,,~~(h_\eta=\lan \eta|F\ran
)\,,
$$
where the higher order terms in ${\mathcal P}$ are omitted.
The highest order of ${\mathcal   P}$
in $\Om$ is called {\em the rank} of the BRST operator $Q$.
If ${\mathcal   A}$ is a Lie algebra defined together with its canonical
action on ${\mathcal   R}$ then $Q$ has the rank one or less. In this case
the BRST operator $Q$ is the extension of the Cartan-Eilenberg operator
giving rise to the cohomology of ${\mathcal A}$ with coefficients in
$\clO({\mathcal R})$. Due to the Jacobi identity
the first two terms in the previous expression provide
the nilpotency of $Q$. It turns out that for the Hamiltonian algebroids ${\mathcal   A}^H$
$\Om$ has the same structure as for  the Lie algebras,
 though the Jacobi identity (\ref{5.5}) has additional terms.
\begin{theor}
The BRST operator $Q$ for the Hamiltonian algebroid ${\mathcal   A}^{H}$
has the rank one:
\beq{5.20}
\Om=\lan \eta|F\ran   +\oh\lan \lb\eta,\eta'\rb|{\mathcal   P}\ran \,.
\end{equation}
\end{theor}
{\sl Proof}.\\
 Straightforward  calculations show that
$$
\{\Om,\Om\}=
\{h_{\eta_1},h_{\eta_2}\}+\oh\lan \lb\eta_2,\eta_2'\rb|F\ran
-\oh\lan \lb\eta_1,\eta_1'\rb|F\ran
$$
$$
+\oh\{h_{\eta_1},\lan \lb\eta_2\eta_2'\rb|{\mathcal   P}_2\ran   \}
-\oh\{h_{\eta_2},\lan \lb\eta_1,\eta_1'\rb|{\mathcal   P}_1\ran   \}+
\f1{4}\{\lan \lb\eta_1,\eta_1'\rb|{\mathcal   P}_1\ran,\lan \lb\eta_2,\eta_2'\rb|{\mathcal
P}_2\ran\}\,.
$$
The sum of the first three terms vanishes due to (\ref{5.12}).
The sum of the rest terms is the left hand side of (\ref{5.5}).
The additional dangerous term may come from the Poisson brackets of
the structure functions
$\{\lb\eta_1,\eta_1'\rb,\lb\eta_2,\eta_2'\rb\}$. In fact, these
brackets vanish because the structure functions
do not depend
on the ghost momenta. Thus, (\ref{5.5}) leads to the desired identity $\{\Om,\Om\}=0$.
$\Box$


\section{Lie algebroids and Poisson sigma-model}
\setcounter{equation}{0}

\subsection{Cotangent bundles to Poisson manifolds
as Lie algebroids}

Let $M$ be a Poisson manifold with the Poisson bivector $\pi=\pi(\ve,\ve')$,
where    $\ve,\ve'$ are sections of the bundle $T^*M$.
It is a skewsymmetric tensor with vanishing Schouten brackets (the Jacobi identity)
$[\pi,\pi]_S=0$.
In local coordinates $x=(x_1,\ldots,x_n)$ it means
\beq{2.1}
\p_i\pi^{jk}(x)\pi^{im}(x)+{\rm c.p.}(j,k,m)=0\,.
\end{equation}
The Poisson brackets are defined on the space
$\clO (M)$
$$
\{f(x),g(x)\}:=\lan dg|\pi| df\ran,~~df,dg\in \G(T^*M)\,.
$$
The Poisson bivector gives rise to  the map
\beq{2.2a}
V^\pi:T^*M\to TM\,,~~
V^\pi_\ve=\lan\ve|\pi|\p\ran\,,~~\ve\in\G(T^*M)\,.
\end{equation}
In particular,
\beq{dx}
\de_\ve x^k= \ve_j\pi^{jk}(x)\,.
\end{equation}
In this way we obtain a map from the space  $\clO(M)$
 to the space of the Poisson vector fields
\beq{2.2}
f\to V_f= \p_if\pi^{ik}\p_k=\lan df|\pi|\p\ran\,,~(\p_i=\frac{\p}{\p x^i})\,.
\end{equation}
The Poisson brackets  can be rewritten as
$\{f(x),g(x)\}=-i_{V_{f}}dg$.

Define brackets on the sections $\ve,\ve'\in\G(T^*M)$
\beq{2.3}
\lb\ve,\ve'\rb=d\lan\ve|\pi(x)|\ve'\ran+\lan d\ve|\pi|\ve'\ran
+\lan\ve|\pi|d\ve'\ran\,,
\end{equation}
$$
\lb\ve,\ve'\rb_k=\ve_j\p_k(\pi^{ji})\ve'_i\,.
$$
\begin{lem}
The  brackets (\ref{2.3}) are  the Lie brackets.
$T^*M$ is a Lie algebroid ${\mathcal   A}$ over the Poisson manifold $M$
with the Lie brackets (\ref{2.3}) and the anchor (\ref{2.2a}).
\end{lem}
{\sl Proof.}\\
It follows from the Jacobi identity (\ref{2.1}), that
the brackets (\ref{2.3}) are the Lie brackets. The
commutator of the vector fields  satisfies (\ref{5.1})
$$
[ V_\ve,V_{\ve'}]=V_{\lb\ve,\ve'\rb}\,.
$$
The property
(\ref{5.2}) follows from the definition of Lie brackets (\ref{2.3}).
$\Box$

\bigskip

The structure functions  of $T^*M$  are defined by the Poisson bivector
$$
f^{jk}_i(x)=\p_i\pi^{jk}(x)\,.
$$
This type of the Lie brackets correspond to a particular choice of the
Dirac structure in the Courant brackets on $TM\oplus T^*M$
\cite{Cou,CouWe}.

\begin{rem}
A linear space $M$ with the linear Poisson brackets $\pi^{jk}(x)=f^{jk}_ix^i$
can be identified with a Lie coalgebra.
 Then  the map $V^\pi$ (\ref{2.2a})
is the coadjoint action
\beq{co}
V^\pi_\ve\sim{\rm ad}^*_\ve \,,
\end{equation}
and the brackets (\ref{2.3}) are the Lie brackets on $M^*$.
\end{rem}


\subsection{Poisson sigma-model and Lie algebroids}

Consider a set $\bfM_{S^1}$ of smooth maps $\tilde{X}\,:\,S^1\to M$
$$
{\bf M}_{S^1} = \{\tilde{X}(t)\,:\,S^1\to M\,,~|t|=1\,,~
\tilde{X}\in C[t,t^{-1}]\otimes\Om^{(m)}(S^1)\,\}\,.
$$
Assume that there exists a Poisson bivector $\pi$ defined on $\bfM_{S^1}$,
with a holomorphic dependence on $\tilde{X}$ and
it is a $(2m-1)$-form on $S^1$
$$
\pi=\pi(\tilde{X})\,,~ \pi\in\wedge^2(T\bfM_{S^1})\otimes \Om^{(2m-1)}(S^1)\,.
$$
By means of $\pi$ define a Lie algebroid $\clA_{S^1}$ over ${\bf M}_{S^1}$, as it was described above.
The sections of the algebroid $\G({\mathcal   A}_{\bf M_{S^1}})=\{\ve_r\}$
are defined by the pairing
\beq{pair}
\oint_{S^1}\de\tilde{\ve}_j\cdot\de \tilde{X}^j\,.
\end{equation}
Therefore, $\tilde{\ve}$ are $(1-m)$-forms on $S^1$ taking values
in the pull-back by $\tilde{X}$ of the cotangent bundle $T^*M$
$$
\G({\mathcal   A}_{\bf M_{S^1}})= \tilde{X}^*(T^*(M))\otimes\Om^{(1-m)}(S^1)\,.
$$

According with (\ref{dx})) and (\ref{2.3}) the anchor and of the brackets  take the form
\beq{7.4c}
\de_\ve \tilde{X}=\lan\ve|\pi(\tilde{X})\,,~~(\de_\ve \tilde{X}^k=\ve_j\pi^{jk}(\tilde{X})\frac{\de}{\de\tilde{X}^k})\,,
\end{equation}
\beq{2.3a}
\lb\tilde{\ve},\tilde{\ve}'\rb=d\lan\tilde{\ve}|\pi(\tilde{X})|\tilde{\ve}'\ran+\lan d\tilde{\ve}|\pi(\tilde{X})|\tilde{\ve}'\ran
+\lan\tilde{\ve}|\pi(\tilde{X})|d\tilde{\ve}'\ran\,.
\end{equation}


Consider a complex curve $\Si_g$ of genus $g$. Let $X(z,\bz)$ and $\ve(z,\bz)$
be smooth continuations of $\tilde{X}(t)$ and $\tilde{\ve}(t)$ from $S^1\subset\Si_g$ on $\Si_g$
and $\bfM$ be the set of the maps
$$
X\,:\,\Si_g\to M\,,~\tilde{X}=X|_{S^1}\,,
$$
$$
{\bf M} = \{X\,:\,C^\infty(\Si_g\to M\,)\,,\\
\otimes\Om^{(m,0)}(\Si_g)\,\}\,.
$$
$$
\ve\in \clG=X^*(T^*(M))\otimes\Om^{(1-m,0)}(\Si_g)\,,~~\tilde{\ve}=\ve|_{S^1}
$$

\bigskip
{\bf The main assumption.}
{\sl The anchor (\ref{7.4c}) and the Lie brackets (\ref{2.3a}) are defined on the
maps $X\in \bfM$ and $\ve\in\clG$. Therefore there exists the Lie algebroid $\clA$ over
$\bfM$ with $\clG$ as a space of sections.}


\subsection{Global Hamiltonian algebroid}

Let   $\xi$ be a $(1-m,1)$-form on $\Si_g$ taking values in sections of the affine space
over the pull-back by
$X$ of the cotangent space  $T^*M$.
$$
\xi\in \G({\bf R}X^*(T^*M))\otimes\Om^{(1-m,1)}(\Si_g)\,.
$$
The affine space $=(\xi,X)$ plays the role of the phase space of the
$2+1$ sigma-model
$$
\Si_g\times \mR\to M\,.
$$
It is endowed  with the canonical symplectic form
\beq{7.2a}
\om=\int_{\Si_g}\lan D X \we D\xi \ran\,.
\end{equation}
Let $c(X,\ve)$ be a cocycle
\beq{7.7}
c(X,\ve)=\int_{\Si_g}\lan\ve|\bp X\ran\,.
\end{equation}
Thereby, one can define the shifted anchor
$$
\de_\ve X^k=\ve_j\pi^{jk}(X)+c(X,\ve)\,.
$$
In principle the cocycle can be trivial and the shift can be removed
(see (\ref{5.11}) and (\ref{5.7})).

The canonical transformations of $\om$ (\ref{7.2a}) according with
 (\ref{5.15a}) and (\ref{5.16a}) are
represented by
\beq{7.40}
\de_\ve X=\lan\ve|\pi(X)\,,
\end{equation}
(see (\ref{7.4c})), and by
\beq{7.5}
\rho_\ve\xi=dc(X,\ve)
+\clL_\ve\xi\,,
\end{equation}
$$
(\rho_\ve\xi_k=\bp\ve_k+\p_k(\ve_n\pi^{nj})\xi_j+\ve_n\pi^{nj})\p_k\xi_j)\,,
$$
where $\clL_\ve$ is the Lie derivative with the vector fields $\lan\ve|\pi(X)$.

The  transformations (\ref{7.40}), (\ref{7.5}) and the brackets (\ref{2.3a}) are
consistent with the orders of forms on $\Si_g$:

\bigskip
\beq{tab}
\begin{tabular}{|c|c|c|c|}
  \hline
  $\ve$ &$ X $& $\xi$ & $\pi $ \\
  \hline
  &    &    &    \\
  $(1-m,0)$ & $(m,0)$ & $(1-m,1)$ & $(2m-1,0)$ \\
  &    &    &    \\
  \hline
\end{tabular}
\eq
\bigskip

Assume that $M$ is flat with global coordinates $(x^j)$, $\,(j=1,\ldots, d)$.
Then another consistent assignment depending on coordinates will be used in Section 5.

\beq{tab1}
\begin{tabular}{|c|c|c|c|}
  \hline
  $\ve_j$ &$ X^j $& $\xi_j$ & $\pi^{jk} $ \\
   \hline
         &    &    &    \\
  $(1-m_j,0)$ & $(m_j,0)$ & $(1-m_j,1)$ & $(m_j+m_k-1,0)$ \\
   &    &    &    \\
  \hline
\end{tabular}
\eq

\bigskip

In Section 6 we will use another pairing between the sections of algebroids
(the covectors) and the tangent to the base vectors (\ref{pa1})
It leads to another correspondence between the forms
and the fields:
\beq{tab2}
\begin{tabular}{|c|c|c|c|}
  \hline
  $\ve$ &$ X $& $\xi$ & $\pi $ \\
  \hline
  &    &    &    \\
  $(-m,0)$ & $(m,0)$ & $(-m,1)$ & $(2m,0)$ \\
  &    &    &    \\
  \hline
\end{tabular}
\eq
\bigskip

The transformations (\ref{7.40}), (\ref{7.5})
 are generated by the first class constraints
 \beq{7.3}
F:=\bp X+\pi(X)|\xi\ran=0\,.
\end{equation}
If we use (\ref{tab}) then $F$ is a $(m,1)$-form on $\Si_g$.
In the case (\ref{tab1}) we come $d=\dim M$ constraints
\beq{7.3a}
F^j=\bp X^j+\pi^{jk}(X)\xi_k=0\,,~~(\dim\,(F^j)=(m_j,1)\,.
\eq

Due to Lemma 2.1 the action (\ref{7.5}) means
the lift of the anchor action from $\bfM$
to  ${\bf  R}/T^*\bfM$ by the cocycle (\ref{7.7}).
The canonical transformations
 are the Hamiltonian transformation
\beq{anc}
\de_{h_{\ve}} f(X,\xi)=\{h_{\ve},f(X,\xi)\}\,.
\end{equation}
Here the Poisson brackets are inverse to the symplectic form $\om$
(\ref{7.2a})
and
\beq{7.8}
h_{\ve}=\int\lan\ve| F\ran=\lan\de_\ve X|\xi\ran+c(X,\ve)\,.
\end{equation}
(see (\ref{5.17})).

Summarizing, we have defined the symplectic manifold
 ${\mathcal   R}=\{(\xi,X)\}$ and the hamiltonian action of
 the algebroid $\clA$ defined by the Hamiltonian (\ref{7.8}).


\subsection{Deformation of complex structure on complex curves}

Following our approach we interpret the constraints  (\ref{7.3}) as consistency
conditions for a linear system. In this and next subsections we take the order
of forms from (\ref{tab}). The passage to  (\ref{tab1}) and  (\ref{tab2}) is straightforward.

Let $\psi\,,\varphi$ be sections of
 $V_m=X^*(T^*M)\otimes\Om^{(-m+1,0)}(\Si_g)$ and \\
$V'_m=X^*(TM)\otimes\Om^{(m,0)}(\Si_g)$, and $B$ is  a linear map $V_m \to V_m'$
\beq{7.4.c}
B(X)\psi=\varphi\,,~~B(X)=\pi (X)\,.
\end{equation}
Define, in addition, two maps
$$
A\,:\, V_m\to V_m\otimes \Om^{(0,1)}(\Si_g)\,,~~
A^*\,:\, V'_m\to V'_m\otimes \Om^{(0,1)}(\Si_g)\,,
$$
\beq{7.5c}
A=-\bp+d\pi(X)|\xi\ran\,,~~A^*=-\bp-d\pi(X)|\xi\ran\,.
\end{equation}
Locally, the operators are defined as
$$
A_k^l\varphi^k=\left(-\bp\de_k^l-\frac{\de}{\de X^k}\pi^{ln}\xi_n\right)\varphi^k\,,
$$
$$
(A^*)_k^l\psi_l=\left(-\bp\de_k^l+\frac{\de}{\de X^k}\pi^{ln}\xi_n\right)\psi_l\,.
$$
Consider the linear system
\beq{7.4a}
\pi (X)\psi=0\,,
\end{equation}
\beq{7.5a}
\left(-\bp-d\pi(X)|\xi\ran\right)\varphi=0\,,
\end{equation}
\beq{7.5b}
\left(-\bp+d\pi(X)|\xi\ran\right)\psi=0\,,
\end{equation}

\begin{lem}
Let the Poisson bivector satisfies the non-degeneracy condition:\\
  $\det\left(d\pi|\psi\ran\right)\neq 0$, $\,(\det\frac{\de}{\de X^i}\pi^{jm}(\psi)_m\neq 0)$
  for any $\psi\in X^*(T^*\bfM)$.
 Then constraints  (\ref{7.3}) are the consistency
conditions for (\ref{7.4a}),(\ref{7.5a}) and (\ref{7.5b}).
\end{lem}
{\sl Proof.}
The consistency condition of these equations is the operator equation
$BA-A^*B=0$ for $B$ (\ref{7.4c}), $A$,   $\,A^*$, (\ref{7.5c}).
After substitution the expressions for $A,A^*,B$ and applying the Jacobi identity
(\ref{2.1}) one comes to the equality
$$
(\bp X^i+\pi^{is}\xi_s)\frac{\de}{\de X^i}\pi^{jm}(\psi)_m=0.
$$
The later is equivalent to the constraint equation (\ref{7.3}) if $\pi$ is
non-degenerate in the above sense. $\Box$

\begin{rem}
The spaces $V_m$ and $V_m'$ are analogs of
coadjoint and adjoint spaces.  In the first example in next Section they
coincide with coadjoint and and adjoint representations of $\sln$.
In what follows we shall consider "vector representations".
\end{rem}

The equations (\ref{7.5a}) and (\ref{7.5b}) define the generalized deformations of the
operator
$\bp$ on $\Si_g$ acting in the space of sections of $V_m$ and $V'_m$.
This deformation is provided by the Poisson bivector $\pi$
and by sections $\xi$ of the affine bundle. We shall apply this scheme for
the concrete Poisson structures below.

\subsection{BRST construction}

Let $G$ be the Lie groupoid corresponding to the Lie algebroid $\clA_\bfM$,
and\\
${\bf R}^{red}={\bf R}//G$ is the corresponding symplectic quotient.

The symplectic quotient ${\bf R}^{red}$ can be described by the BRST technique.
Define the BRST anticommuting ghosts
$$
\eta\in\G(X^*(T^*M)\otimes\Om^{(1-m,0)}(\Si_g)\,,
$$
and their momenta
$$
{\mathcal   P}\in\G(X^*(TM)\otimes\Om^{(m,1)}(\Si_g)\,.
$$
The classical BRST complex is the set of fields
\beq{BRS}
\bigwedge{}^\bullet\left
(\G(X^*(TM)\otimes\Om^{(1-m,0)}(\Si_g))
\oplus\G(X^*(T^*M)\otimes\Om^{(m,1)}(\Si_g))
\right)\otimes \clO({\mathcal   R})\,.
\end{equation}

Theorem 2.1 states that the BRST operator has the rank one
\beq{OM}
\Om=\int_{\Si_g}\lan\eta |F\ran+
\int_{\Si_g}\lan\lb\eta,\eta\rb|{\mathcal   P}\ran\,.
\end{equation}

Remind that the classical observables on ${\bf R}^{red}$ are elements from $H^0(Q)$.
 The moduli space of deformations of complex structures is a part of
${\mathcal   R}^{red}$.

\bigskip
Let us briefly repeat the steps that lead to the moduli space of
deformations
of complex structure on the disk $\Si_g$.
\begin{itemize}
  \item We start with a Poisson manifold $M$ and define the Lie algebroid ${\mathcal A}$
  over $M$ (Lemma 3.1).

  \item This algebroid has infinite-dimensional version ${\mathcal A}_{\bfM}$ if one
consider the maps $X$
  from $\Si_g$ to $M$.

  \item We define the set of fields $(X,\xi)$, where $\xi\in X^*(T^*M)$.
  It is the affine space ${\bf R}/T^*{\bf M}$.

  \item The anchor action and the representation of ${\mathcal   A}$ in $\bf R$ is generated
  by the first class constraints (\ref{7.3}).
  They are the compatibility conditions for the linear system (\ref{7.4a}), (\ref{7.5a}) and
(\ref{7.5b}).
  Two last equations define the generalized deformation of complex structures on
in the space of sections  $V_m$, $V'_m$.

  \item The reduced phase space ${\bf R}^{red}$ can be described in terms of the BRST complex (\ref{BRS}) with  $\Omega$ (\ref{OM}).
\end{itemize}

We repeat these steps in concrete cases considered in the rest part of paper.


\section{Two examples of Hamiltonian algebroids with
Lie algebra symmetries}
\setcounter{equation}{0}

In this section we consider two examples, where the spaces of sections
of the gauge algebras are replaced by sections of Lie algebroids,
though the results can be obtained within standards Lie algebras symmetries
implying a trivial anchor action.
We  construct Lie algebroids to illustrate our approach.

Let $\Si_{g,n}$ be a complex curve of
 genus $g$
with $n$ marked points. The first example  is the moduli
space of flat bundles over $\Si_{g,n}$. It will become clear later, that it is
an universal system containing hidden algebroid symmetries.
 The second example is the moduli space of
the projective structures ($ \clW_2$-structures) on $\Si_{g,n}$.

The generalization of the latter example is the $\clW_N$-structures, where
the symmetries are
defined by a nontrivial Lie  algebroid, will be considered
in last Sections.


\subsection{Flat bundles  with regular singularities}

We consider a rank $N$ trivial vector bundle  $E$ over
 $\Si_{g,n}$. Define the derivatives
 $d': E\to E\otimes\Om^{(1,0)}(\Si_{g,n})$,
$d'': E\to E\otimes\Om^{(0,1)}(\Si_{g,n})$.


\subsubsection{Local Lie algebroid}

On a disk  $D\subset\Si_{g,n}$ one can choose the derivatives
 in the form
\beq{8.0}
d'=(\ka\p +A)\otimes dz\,, ~~d''=\bp+\bA\,,~~\ka\in\mC\,
\end{equation}
where $\p=\p_z$, $\bp=\p_{\bz}$, $z$ is a local coordinate and $A(z,\bz)\,$,
$\,\bA(z,\bz)$ are $\sln$ valued $C^\infty(D)$ functions.
Let  $M_{SL_N}(S^1)$  be  the set  $\{d'\}$ restricted on the boundary
$S^1$ of $D$.
It has a structure of the affine Lie coalgebra
$\hat{L}^*(\sln)$ with the Lie-Poisson brackets on the space of functionals
$\clO(M_{SL_N}(S^1))$
$$
\{f(A),g(A)\}=
\oint_{S^1}\tr\left([df(A),dg(A)]A)+df(A)\p (dg(A)) \right)\,,~~\{f(A),\ka\}=0\,,
$$
where $df(A)$ is a variation of $f(A)$.
Thereby, $M_{SL_N}(S^1)$ can be considered as
 the base of the Lie algebroid ${\mathcal   A}_{SL_N}(S^1)$
 (see Remark 3.1). It corresponds to $m=1$ in (\ref{tab}).

Following (\ref{2.3}) we define
the space of sections of the algebroid. It is the  Lie algebra of smooth functionals
${\mathcal   G}_{SL_N}(S^1)=\hat{L}(\sln)$ with coefficients from $\clO(M_{SL_N}(S^1))$.
The variable dual to $\ka$ corresponds to the central charge of $\hat{L}(\sln)$.
In this way we come to  the Lie brackets
\beq{lbra}
[(\ve(t),0),(\ve'(t')),0)]_{jk}=(\de(t,t')\sum_{l=1}^N\ve_{jl}(t)\ve'_{lk}(t')-\ve'_{jl}(t')
\ve_{lk}(t), \oint_{S^1}\tr(\ve(t)\p\ve'(t)))\,,
\end{equation}
where $\de(t,t')=\sum_{i\in\mZ}t^i(t')^{-i-1}dt$.
The integral defining the central extension represents a nontrivial
2-cocycle (see (\ref{5.8})).

The anchor is defined by the gauge transformations
\beq{8.1}
\de_\ve A=\ka\p\ve+[A,\ve]\,.
\end{equation}
Note, that the multiplication of sections on functionals from $\clO(M_{SL_N}(S^1))$
modifies the brackets (\ref{lbra}) according with (\ref{5.2})
Similarly, we have for the gauge transformations
$$
\de_{f(A)\ve} A=f(A)(\ka\p\ve+[A,\ve])\,.
$$


\subsubsection{Global Lie algebroid}

Though the Poisson structure is defined only on $S^1$,
the gauge algebra (\ref{lbra}) and the gauge transformations (\ref{8.1}) are well
defined on the curve. Thereby it is possible to define the global Lie algebroid.
We specify the behavior of the fields in neighborhoods
of the marked points.
Assume that $A$ has first order holomorphic poles at the marked points
\beq{8.0a}
A|_{z\to x_a}=\frac{A_a}{z-x_a}\,.
\end{equation}
We denote by ${\mathcal   G}_{SL_N}=\{\ve\}$  the Lie algebra
of the smooth gauge transformations on  $\Si_{g,n}$. Assume that at the marked points
the gauge transformations are nontrivial
\beq{8.0b}
\ve|_{z\to x_a}=r_a+O(z-x_a)\,, ~~r_a\in\sln\neq 0\,.
\end{equation}

At the marked points we add a collection $P$ of $n$  elements from  the
Lie coalgebra ${\rm sl}^*(N,{\mC})$
$$
P=\left\{\bfp=(p_1,\ldots,p_n)\right\}\,,
$$
endowed with the Lie-Poisson structure
$
\{f(p_a),g(p_b)\}=\de_{ab}\lan p_a|[ df,dg]\ran$.

The gauge algebra ${\mathcal   G}_{SL_N}$ acts on $P$ by the evaluation maps as
\beq{8.1a}
\de_\ve p_a=[p_a,r_a]\,,~~
\ve\in {\mathcal   G}_{SL_N}\,.
\end{equation}
Thereby, we have defined a Lie algebroid
${\mathcal   A}_{SL_N} ={\mathcal   G}_{SL_N}\times M_{SL_N}$
over\\
$M_{SL_N}=\{d',P\}$  with the anchor map (\ref{8.1}), (\ref{8.1a}).

The cohomology
$H^i({\mathcal   A}_{SL_N})=H^i({\mathcal   G}_{SL_N},\clO(M_{SL_N}))$
are the standard cohomology
of the gauge algebra ${\mathcal   G}_{SL_N}$ with the cochains taking values
in holomorphic functionals on $M_{SL_N}$.  There is a nontrivial one-cocycle
\beq{8.2}
c(A,\bfp;\ve)=\int_{\Si_{g,n}}\tr\,\left(\ve(
\bp A-2\pi i\sum_{a=1}^n\de(x_a)p_a)
\right)
\end{equation}
$$
=
\lan\ve|\bp A\ran-2\pi i\sum_{a=1}^n\tr(r_a\cdot p_a)
$$
representing an element of $H^1({\mathcal   A}_{SL_N})$.
 This cocycle provides a nontrivial
extension of the anchor action (see (\ref{5.11}))
$$
\hat{\de}_\ve f(A,\bfp)
=\lan\ve|\bp A-\p (df(A))+[df(A),A]\ran
-2\pi i\sum_{a=1}^n\tr(r_a p_a) \,.
$$

Next consider $2g$ contours $\ga_\al,~(\al=1,\ldots,2g)$
 generating $\pi_1(\Si_{g})$. The contours determine the 2-cocycles
\beq{8.3}
c_\al(\ve_1,\ve_2)=\int_{\ga_\al}\tr(\ve_1\p\ve_2)
\end{equation}
(see (\ref{lbra}).
The cocycles (\ref{8.3}) lead to $2g$ central extensions
$\hat{\mathcal   G}_{SL_N}$ of ${\mathcal   G}_{SL_N}$
$$
\hat{\mathcal   G}_{SL_N}={\mathcal   G}_{SL_N}\oplus_{\al=1}^{2g}{\bf C}\La_\al\,,
$$
$$
[(\ve_1,0),(\ve_2,0)]_{c.e.}=
\left([\ve_1,\ve_2],\sum_{\al}c_\al(\ve_1,\ve_2)\right)\,.
$$


\subsection{Global Hamiltonian algebroid}

To define the corresponding Hamiltonian algebroid we
consider the cotangent bundle $T^*M_{SL_N}(\Si_{g,n})$ with the sections
$\bar{\Phi}\in\Om^{(0,1)}(\Si_{g,n},{\rm sl}(N,{\mC}))$.
\begin{rem}
The form $\bar{\Phi}$ is dual to $d'$.
Here and in what follows we do not consider the dual to $\ka$ variables.
\end{rem}
Define
the affine space ${\mathcal   R}^0_{SL_N}$  over $T^*M_{SL_N}(\Si_{g,n})$
as the space of sections   $\{d''=\bp+\bA\}$.

   The symplectic form on
${\mathcal   R}_{SL_N}$ is
\beq{om0}
\om^0= \int_{\Si_{g,n}}\tr( DA\wedge D\bA)=\lan DA\wedge D\bA\ran\,.
\end{equation}

Consider the contributions of the marked points in the symplectic structure.
 We define there the symplectic manifold
$$
(T^*G_1\times\ldots,\times T^*G_n)\,.
$$
with the form
\beq{8.4c}
\sum_{a=1}^n\om_a=\sum_{a=1}^n\lan(D(p_a g_a^{-1})\wedge Dg_a\ran\,.
\end{equation}
Here $\om_a$ is the canonical symplectic form on $T^*G_a\sim T^*\SLN$.
We pass from $T^*G_a$ to the coadjoint orbits
\beq{orb}
{\mathcal   O}_a=\{p_a=g_a^{-1}p_a^{(0)}g_a~|~p_a^{(0)}=\di(\la_{a,1},\ldots,\la_{a,N})\,,~
~g_a\in{\rm SL}(N,{\bf C})\}\,.
\end{equation}
and assume that the orbits are generic $\la_{a_j}\neq\la_{a_k}$, for $j\neq k$.
The orbits are the symplectic quotient
${\mathcal   O}_a\sim \SLN\backslash\backslash  T^*G_a$ with respect to the action
$g_a\to f_ag_a,~f_a\in\SLN$.
The form $\om_a$ coincides on ${\mathcal   O}_a$ with the Kirillov-Kostant form
$\om_a=\lan D(g_a^{-1} p^{(0)}_a)\wedge Dg_a\ran$. It was mentioned above (Example 1,\, in {\bf 2.4})
that the orbits ${\mathcal   O}_a$ are the affine spaces $Aff(T^*Fl_a(N))$ over
the cotangent bundles $T^*Fl_a(N)$ to the flag varieties $Fl_a(N)$.

Eventually we come to the symplectic manifold
$$
{\mathcal R}_{SL_N}=
({\mathcal R}^0_{SL_N};{\mathcal O}_1\times\ldots{\mathcal O}_n)
\sim (Aff(T^*E);Aff(T^*Fl_1)\times\ldots Aff(T^*Fl_n))\,,
$$
\beq{8.4a}
\om=\om^0+\sum_{a=1}^n\om_a=\int_{\Si_g}\tr( DA\wedge D\bA)+\sum_{a=1}^n\tr(D(g_a^{-1}p^0_a)\wedge
Dg_a )\,.
\end{equation}
According with (\ref{7.7}) the pass from $T^*M_{SL_N}(D^\times)$
to ${\mathcal   R}_{SL_N}$ is provided
by the cocycle (\ref{8.2})



 Consider the Hamiltonian
$$
h_\ve=
\int_{\Si_g}\tr\,\ve\left( F(A,\bA)-2\pi i\sum_{a=1}^n\de(x_a)p_a)\right),~~
F(A,\bA)=\bp A-\ka\p\bA+[\bA,A]\,.
$$
The Hamiltonian generates the canonical vector fields (\ref{8.1}) and
\beq{cantr}
\rho_\ve\bA=\bp\ve+[\bA,\ve]\,,~~\rho_\ve g_a=g_ar_a\,.
\end{equation}
(see (\ref{5.16a})).
The global version of this transformations is the gauge groupoid
$G_{SL_N}$ acting on ${\mathcal   R}_{SL_N}$.
The flatness condition
\beq{8.5a}
F(A,\bA)-2\pi i\sum_{a=1}^n\de(x_a)p_a=0
\end{equation}
is the moment constraint with respect to this action.

The flatness is the compatibility condition for the linear
system
\beq{8.5c}
\left\{
\begin{array}{c}
(\ka\p +A)\psi=0\,,\\
(\bp +\bA)\psi=0\,,
\end{array}
\right.
\end{equation}
where $\psi\in\Om^0(\Si_{g,n},{\rm Aut} E)$ as in {\bf Lemma\, 3.2}.
We can consider the same system for  $\psi\in\Om^0(\Si_{g,n}, E)$
(see {\bf Remark\,3.1}).

The second equation
describes the deformation of the holomorphic structure on $E$.

\subsubsection{The moduli space of flat bundles}

The moduli space ${\mathcal   M}^{flat}_N$ of flat $\SLN$-bundles
is the symplectic quotient
${\mathcal   R}_{SL_N}// G_{SL_N}$. It has  dimension
\beq{8.5b}
\dim {\mathcal   M}^{flat}_N=2(N^2-1)(g-1)+N(N-1)n\,,
\end{equation}
where the last term is the contribution of the coadjoint orbits ${\mathcal   O}_a$.

As in the general case $\clO({\mathcal   M}^{flat}_N)$ can be
identified with cohomology group $H^0(Q)$ of the BRST operator $Q$ which we are going to define.
Let $\eta\in\Om^{(0)}(\Si_{g,n},{\rm End}E)$ be the ghost field
and ${\mathcal   P}$ is its momentum ${\mathcal   P}\in\Om^{(1,1)}(\Si_{g,n},{\rm End}E)$.
 Consider the algebra
$$
\clO({\mathcal   R}_{N})\otimes\wedge^\bullet
 \left({\mathcal   G}_{ SL_N})\oplus {\mathcal   G}^*_{ SL_N}\right)\,.
$$
Then the BRST operator $Q$ acts on functionals on this algebra as
$$
Q\Psi(A,\bA,p_a,\eta,{\mathcal   P})=\{\Om,\Psi(A,\bA,p_a\eta,{\mathcal   P})\}\,,
$$
 where
$$
\Om=\int_{\Si_g}\tr\,\eta \left(F(A,\bA)-2\pi \imath\sum_{a=1}^np_a\de(x_a)\right)\ran+
\oh\int_{\Si_g}\tr( [\eta,\eta']|{\mathcal
P})   \,.
$$


\subsection{Projective structures on $\Si_{g,n}$}

\subsubsection{Local Lie algebroid}

The set of
projective connections  $M_2(D)$ on a  disk $D\subset\Si_{g,n}$
 is  represented by the second order differential operators $\ka^2\p^2-T$,
 where $T=T(z,\bz)\in C^\infty(D)$ $\,T(z,\bz)\in\Om^{(2,0)}(D)\,$, $\,\ka\in\mC$.

  The set $M_2(S^1)$  $(S^1\sim\p D)$   is a Poisson manifold
with the brackets
\beq{vir}
\{T(t),T(s)\}=\left(-\frac{1}{2}\ka^3\p^3+2T\ka\p+\ka\p T\right)\de(t,s)\,,~~ \{T,\ka\}=0\,,
\end{equation}
where $\p=\p_t$ and $\de(t,s)=\sum_{k\in\mZ}t^ks^{-k-1}dt$.
This case corresponds to (\ref{tab}) $m=2$

 The dual space with respect to the pairing
 $$
 \lan(T,\ka)|(\ve,c)\ran=\oint_{S^1}\de\ve\cdot \de T+\ka c
 $$
 is the central extended Lie algebra of vector fields ${\mathcal G}_2(S^1)=\{(\ve,c)\}$ on $S^1$
 $$
(\ve,c)\,,~~ \ve=\ve(z)\frac{\p}{\p z}\in\Om^{(-1,0)}(S^1)\,,~~
\ve(z)\in  C^\infty(S^1)\,.
$$
   The commutation relations can be read off
 from the Poisson brackets  (see (\ref{2.3}))
\beq{8.7b}
[(\ve_1,0),(\ve_2,0)]=
(\ve_1\kappa\p\ve_2-\ve_2\kappa\p\ve_1,\oint_{S^1}\ve_1(-\frac{3}{2}\ka^2\p^3+2T\p+\p T)\ve_2)\,.
\end{equation}
One can omit the last two terms under the integral since they form an exact cocycle
(see (\ref{5.8}) and we deal with the standard cocycle.

The coadjoint action of ${\mathcal G}_2(S^1)$ on $M_2(S^1)$
\beq{8.7}
\de_{\ve}T(z,\bz)=-\ve\ka\p T-2T\ka\p\ve+\frac{1}{2}\ka^3\p^3\ve
\end{equation}
defines the anchor in the vector bundle ${\mathcal   A}_{2}(S^1)$
over  $M_2(S^1)$.


\subsubsection{Global Lie algebroid}

The algebroid (\ref{8.7b}), (\ref{8.7}) can be  defined
globally over the space $M_2=M_2(\Si_{g,n})$ of projective connections $T$
on $\Si_{g,n}$. We assume that $T(z,\bz)$ is smooth on $\Si_{g,n}$ and  has
poles at the marked points $x_a,(a=1,\ldots,n)$
  up to the second order:
\beq{8.8a}
T|_{z\rar x_a}\sim\frac{T^a_{-2}}{(z-x_a)^2}+\frac{T^a_{-1}}{(z-x_a)}+\ldots\,,
\end{equation}

The section of the algebroid are smooth chiral vector fields
$$
{\mathcal   G}_2(\Si_{g,n})=\G(\Om^{(-1,0)}(\Si_{g,n}))=\{\ve(z,\bz)\frac{\p}{\p z}\}\,.
$$
Assume that the vector fields have the first order holomorphic nulls at the marked points
\beq{8.8b}
\ve|_{z\rar x_a}=r_a(z-x_a)+o(z-x_a),~~r_a\neq 0\,.
\end{equation}
We denote this global algebroid
${\mathcal A}_{2}={\mathcal G}_2(\Si_{g,n})\oplus M_2(\Si_{g,n})$.

Consider the cohomology $H^\bullet ({\mathcal A}_2)\sim H^\bullet ({\mathcal G}_2,M_2)$.
Due to (\ref{8.7}) and (\ref{8.8b})  $\de_{\ve}T^a_{-2}=0$ and thereby
 $T^a_{-2}$ in (\ref{8.8a}) represents an element
from $H^0({\mathcal   A}_2)$.

The anchor action (\ref{8.7}) can be extended by the one-cocycle $c(T;\ve)$
representing a nontrivial element of $H^1({\mathcal   A}_2)$
\beq{8.9a}
c(T;\ve)=\int_{\Si_{g,n}}\ve\bp T=\lan\ve|\bp T\ran\,,
\end{equation}
$$
\hat{\de}_{\ve}f(T)=\lan\de_\ve T|df(T)\ran+c(T;\ve)\,.
$$
The contribution of the marked points in (\ref{8.9a}) is $2\pi ir_aT_{-2}^a$.

There exist  $2g$ nontrivial two-cocycles defined by the integrals
over non contractible contours $\ga_\al$:
$$
c_\al (\ve_1,\ve_2)=\ka^3\oint_{\ga_\al}\ve_1\p^3\ve_2\,.
$$
The cocycles give rise to the central extension $\hat{\mathcal G}_2$ of the Lie algebra
of the vector fields on $\Si_{g,n}$ (see (\ref{8.7b})).


\subsubsection{Global Hamiltonian algebroid}

The affine space ${\mathcal   R}_2(\Si_{g,n})/T^*M_2(\Si_{g,n})$ has
the Darboux coordinates
$T$ and $\mu$,
where $\mu\in\Om^{(-1,1)}(\Si_{g,n})$ is the Beltrami differential.
The anchor (\ref{8.7}) is lifted to ${\mathcal   R}_2(\Si_{g,n})$ as
\beq{8.11}
\de_\ve\mu=-\ve\ka\p\mu +\mu\ka\p\ve+\bp\ve\,,
\end{equation}
where the last term occurs due to the cocycle (\ref{8.9a}).
The symplectic form on ${\mathcal   R}_2(\Si_{g,n})$ is
\beq{sympl}
\om=\int_{\Si_{g,n}}DT\wedge D\mu\,.
\end{equation}
\begin{rem}
The space ${\mathcal   R}_2$ is the classical phase space
of the $2+1$-gravity on $\Si_{g,n}\times I$ \cite{Ca}. The Beltrami differential
$\mu$ is related to the conformal class of metrics on $\Si_{g,n}$
and plays the role of a coordinate, while $T$ is a momentum.
In our construction the role of $\mu$ and $T$ is interchanged.
\end{rem}

We specify the dependence of $\mu$ on the positions of
the marked points in the following  way. Let ${\mathcal   U}'_a$ be
 neighborhoods
 of the marked points $x_a\,$,\\$(a=1,\ldots,n)$
such that ${\mathcal   U}'_a\cap{\mathcal   U}'_b=\emptyset$ for $a\neq b$.
Define a smooth function $\chi_a(z,\bz)$
\beq{cf}
\chi_a(z,\bz)=\left\{
\begin{array}{cl}
1,&\mbox{$z\in{\mathcal   U}_a$ },~{\mathcal   U}'_a\supset{\mathcal   U}_a\\
0,&\mbox{$z\in\Si_g\setminus {\mathcal   U}'_a$}\,.
\end{array}
\right.
\end{equation}
Due to (\ref{8.11}) at the neighborhoods of the marked points $\mu$ is defined
 up to the term $\bp(z-x_a)\chi(z,\bz)$.
Then $\mu$ can be represented as
\beq{mu}
 \mu=\sum_{a=1}^n[t_{0,a}+t_{1,a}(z-x_a)+\ldots]\mu^0_a\,,~~
\mu^0_a=\bp \chi_a(z,\bz)\,,~~(t_{0,a}=x_a-x_a^0)\,,
\end{equation}
where only $t_{0,a}$ can not be removed by the gauge transformations (\ref{8.8b}),
(\ref{8.11}).

Contribution of the marked points to the symplectic form (\ref{sympl}) takes the form
\beq{8.14}
\sum_{a=1}^nDT^a_{-2}\wedge Dt_{1,a}+DT^a_{-1}\wedge Dt_{0,a}\,.
\end{equation}

The canonical transformations are generated by the Hamiltonian
\beq{8.12}
h_\ve=\int_{\Si_{g,n}}\ve F(T,\mu)=\int_{\Si_{g,n}}\mu\de_\ve T+c(T,\ve)\,,
\end{equation}
where
\beq{8.13}
F(T,\mu)=
(\bp+\mu\ka\p+2\ka\p\mu)T-\frac{1}{2}\ka^3\p^3\mu\,.
\end{equation}
We put $F(T,\mu)=0$
\beq{8.15}
(\bp+\mu\ka\p+2\ka\p\mu)T-\frac{1}{2}\ka^3\p^3\mu=0\,.
\end{equation}

Let $\psi$ be a $(-\oh,0)$ differential. Then (\ref{8.15})
is the compatibility condition for the linear system
\beq{8.16}
\left\{
\begin{array}{l}
(\ka^2\p^2-T)\psi=0\,,\\
(\bp+\mu\ka\p -\oh\ka\p\mu) \psi=0\,.
\end{array}
\right.
\end{equation}
It is analog of the vector representation mentioned in {\bf Remark\,3.1}.
The system (\ref{7.4a}) - (\ref{7.5b}) in this case has the form
$$
\left\{
\begin{array}{l}
(-\oh\ka^3\p^3+2T\ka\p +\ka\p T))\phi=0\,,\\
(\bp+\mu\ka\p -\ka\p\mu) \phi=0\,.
\end{array}
\right.
$$
Here $\phi$ is a section of the adjoint representation $\Om^{(-1,0)}(\Si_{g,n})$.

It follows from the second equations in both systems that the Beltrami differential $\mu$ provides
 the deformation of complex structure on $\Si_{g,n}$.


\subsubsection{The moduli space ${\mathcal   W}_2$}

Let $G_2$ be the group  corresponding
to the Lie algebra ${\mathcal   G}_2$.
\begin{defi}
The moduli space ${\mathcal   W}_2$ of $W_2$-gravity on $\Si_{g,n}$ is
the symplectic quotient of ${\mathcal   R}_2$ with respect to the action of $G_2$,
$$
{\mathcal   W}_2={\mathcal   R}_2//G_2=\{F(T,\mu)=0\}/G_2\,.
$$
\end{defi}
It has dimension $6(g-1)+2n$.
The space of observables is isomorphic
to the cohomology $H^0$ of the BRST complex. It is generated by
the fields $T,\mu\in {\mathcal   R}_2$, the ghosts fields $\eta\in\Om^{(-1,0)}(\Si_{g,n})$
 and the ghosts momenta ${\mathcal   P}\in\Om^{(2,1)}(\Si_{g,n})$.
The BRST operator $Q$ is defined by $\Om$
$$
\Om=\int_{\Si_{g,n}}\eta F(T,\mu)+\oh\int_{\Si_{g,n}}[\eta,\eta']{\mathcal   P}.
$$
The first term is just the Hamiltonian (\ref{8.12}), where the vector fields are
replaced by the ghosts.


\section{Hamiltonian algebroid structure in ${\mathcal   W}_3$-gravity}
\setcounter{equation}{0}

Now consider the concrete example of the general construction with a nontrivial
algebroid structure.
It is the $W_N$ structures on $\Si_{g,n}$ \cite{P,BFK,GLM}.
They generalize the $W_2$ structure described in previous Section.
In this Section we consider in details the $W_3$ case.


\subsection{$\SLN$-opers}

Opers are $G$-bundles over complex curves with additional structures
\cite{Tel,BD}.
Let $E_N$ be a  $\SLN$-bundle over $\Si_{g,n}$. It is a $\SLN$-{\em oper} if
 there exists a flag filtration
$E_N\supset  \ldots\supset   E_1\supset    E_0=0$ and a covariant derivative,
 that acts as
$\nabla:~E_j\subset E_{j+1}\otimes\Om^{(1,0)}(\Si_{g,n})$.
Moreover, $\nabla$ induces an isomorphism
$E_j/E_{j-1}\to E_{j+1}/E_{j}\otimes\Om^{(1,0)}(\Si_{g,n})$.
It means that locally
\beq{opN}
\nabla=\ka\p-
\left(
\begin{array}{cccccc}
0   & 1 & 0 & \ldots&   & 0 \\
0   & 0 & 1 & \ldots&   &   \\
    &   &   & \cdot &   &   \\
 0  & 0 &   &       & 0&  1  \\
W_N& W_{N-1}&\ldots&&W_2&0
\end{array}
\right)\,,
\end{equation}
where the matrix elements $W_k=W_k(z,\bz)$ are smooth.
In other words, we define the $N$-order differential operator on $\Si_{g,n}$
\beq{dopN}
L_N=\ka^N\p^N-W_2\ka^{N-2}\p^{N-2}\ldots-W_N~:
~\Om^{(-\frac{N-1}{2},0)}(\Si_{g,n})\to\Om^{(\frac{N+1}{2},0)}(\Si_{g,n})
\end{equation}
with vanishing subprinciple symbol. The $\GLN$-opers come from the $\GLN$-bundles
and have the additional term $-W_1\p^{N-1}$ in (\ref{dopN}).

In this section we consider $\SLt$-opers and postpone the general case to next Section.
It is possible to choose $E_1=\Om^{-1,0}(\Si_{g,n})$.
For $N=3$ we have
\beq{6.1}
\nabla=\ka\p -\thmat{0}{1}{0}{0}{0}{1}{W}{T}{0}\,,
\end{equation}
and the third order differential operator
\beq{6.2}
L_3=\ka^3\p^3-T\ka\p-W:~\Om^{(-1,0)}(\Si_{g,n})\to\Om^{(2,0)}(\Si_{g,n})\,.
\end{equation}


\subsection{Local Lie algebroid over SL$(3,\mC)$-opers}

Consider the set $M_3(D)=\{L_3\}$ of $\SLt$-opers on a disk $D\subset\Si_{g,n}$.
On $S^1=\p D$
this set becomes a Poisson manifold with respect to the AGD brackets \cite{Ad,GD}
\beq{6.9a}
\{T(t),T(t')\}
=\left(-2\ka^3\p^3+2T(t)\ka\p+\ka\p T(t)\right)\de(t-t')\,,
\end{equation}
\beq{6.10a}
\{T(t),W(t')\}=\left(\ka^4\p^4-T(t)\ka^2\p^2+3W(t)\ka\p-\ka\p W(t)\right)\de(t-t')\,,
\end{equation}
\beq{6.12a}
\{W(t),W(t')\}=
\end{equation}
$$
+\left(\frac{2}{3}\ka^5\p^5-\frac{4}{3}T(t)\ka^3\p^3-
2\p T(t)\ka^2\p^2+\left(\frac{2}{3}T(t)^2-2\ka^2\p^2T(t)+2\ka\p W(t)\right)\ka\p
\right.
$$
$$
\left. +\left(\ka^2\p^2W(t)-\frac{2}{3}\ka^3\p^3T(t)+\frac{2}{3}T(t)\ka\p T(t)\right)
\right)\de(t-t')\,.
$$
\begin{rem}
These brackets can be obtained in two ways. It was pointed in
In Ref.\,\cite{DS} they   derived via the Poisson reduction
from the Lie-Poisson brackets on the Lie coalgebra of the Borel subalgebra $\hat{L}($sl$(3,mC))$ by the
with respect the action of the unipotent subgroup. Another scheme was proposed in \cite{BFK},
where the original coalgebra is $\hat{L}^*($sl$(3,mC))$ and the action is generated by a maximal
parabolic subgroup.
\end{rem}

In this way $M_3(S^1)$ can be considered as
  a base of a Lie algebroid
${\mathcal   A}_3(S^1)\sim  T^*M_3 (S^1)$.
This situation corresponds to (\ref{tab1}) with $M\sim\mC^2$, $m_1=2$
and $m_2=3$.
To define the space of its sections we consider
the dual space  $M^*_3 (S^1)$  of second order
differential operators on $S^1$ with a central extension
\beq{m3}
M^*_3 (S^1)=
\{ (\ve^{(1)}\frac{d}{dt}+ \ve^{(2)}\frac{d^2}{dt^2},c)\}\,.
\end{equation}
This space is defined by the pairing (see (\ref{pair}))
\beq{pa}
\lan(T,W,\ka)|(\ve^{(1)},\ve^{(2)},c)\ran=
\oint_{S^1}(\ve^{(1)}T+ \ve^{(2)}W)+\ka c\,.
\end{equation}
 Following  (\ref{2.3}) we define
the Lie brackets on $ T^*M_3 (S^1)$ by means of the AGD Poisson structure
\beq{6.6}
\lb(\ve^{(1)}_10),(\ve^{(1)}_2,0)\rb=
(\ka(\ve^{(1)}_1\p\ve^{(1)}_2-\ve^{(1)}_2\p\ve^{(1)}_1) \frac{d}{dt}\,,\,c(\ve^{(1)}_1,\ve^{(1)}_2))\,,
\end{equation}
\beq{6.7}
\lb(\ve^{(1)},0),(\ve^{(2),0})\rb=\left(
-\ve^{(2)}\ka^2\p^2\ve^{(1)}) \frac{d}{dt}+
(-2\ve^{(2)}\ka\p\ve^{(1)}+\ve^{(1)}\ka\p\ve^{(2)}) \frac{d^2}{dt^2}\,,\,c(\ve^{(1)}_1,\ve^{(2)}_2)\right)\,,
\end{equation}
\beq{6.8}
\lb\ve^{(2)}_1,\ve^{(2)}_2\rb=\left((
\frac{2}{3}[\ka\p(\ka^2\p^2-T)\ve^{(2)}_1]\ve^{(2)}_2 -
 \frac{2}{3}[\ka\p(\ka^2\p^2-T)\ve^{(2)}_2]\ve^{(2)}_1) \frac{d}{dt} +
 \right.
\end{equation}
$$
\left.
(\ve^{(2)}_2\ka^2\p^2\ve^{(2)}_1-\ve^{(2)}_1\ka^2\p^2\ve^{(2)}_2)\frac{d^2}{dt^2}\,,
 \,c(\ve^{(2)}_1,\ve^{(2)}_2))  \right)\,.
$$
Here $c(\ve^{(j)}_1,\ve^{(k)}_2)$ are the cocycles
\beq{coc1}
c(\ve^{(j)}_1,\ve^{(k)}_2)=\oint_{S^1}\la(\ve^{(j)}_1,\ve^{(k)}_2),
~~(j,k=1,2)\,,
\end{equation}
$$
\la(\ve^{(1)}_1,\ve^{(1)}_2)=-2\ve_1^{(1)}\ka^2\p^3\ve_2^{(1)}+\ldots\,,~~
\la(\ve^{(1)}_1,\ve^{(2)}_2)=\ve^{(1)}_1\ka^3\p^4\ve^{(2)}_2+\ldots\,,
$$
$$
\la(\ve^{(2)}_1,\ve^{(2)}_2)=\frac{2}{3}
\ve_1^{(2)}\ka^4\p^5\ve_2^{(2)}+\ldots\,,
$$
and ellipses means the terms depending on lesser degrees of $\ka$.
It can be proved that $sc=0$
(\ref{5.9}) and that $c$ is not exact.
Note, that the brackets (\ref{6.6}) define the algebra  of the vector fields and  the commutation relations are their generalization to the second order
differential operators.

According with (\ref{dx}) the anchor action in ${\mathcal   A}_3(S^1)$ has the form
\beq{6.9}
\de_{\ve^{(1)}}T=-2\ka^3\p^3\ve^{(1)}+2T\ka\p\ve^{(1)}+\ka\p T\ve^{(1)}\,,
\end{equation}
\beq{6.10}
\de_{\ve^{(1)}}W=-\ka^4\p^4\ve^{(1)}+3W\ka\p\ve^{(1)}+\ka\p W\ve^{(1)}+T\ka^2\p^2\ve^{(1)}\,,
\end{equation}
\beq{6.11}
\de_{\ve^{(2)}}T=\ka^4\p^4\ve^{(2)}-T\ka^2\p^2\ve^{(2)}+(3W-2\ka\p T)\ka\p\ve^{(2)}+
(2\ka\p W-\ka^2\p^2T)\ve^{(2)}\,,
\end{equation}
\beq{6.12}
\de_{\ve^{(2)}}W=\frac{2}{3}\ka^5\p^5\ve^{(2)}-\frac{4}{3}T\ka^3\p^3\ve^{(2)}-
2\ka^3\p T\p^2\ve^{(2)}+
\end{equation}
$$
\ka(\frac{2}{3}T^2-2\ka^2\p^2T+2\ka\p W)\p\ve^{(2)}+
(\ka^2\p^2W-\frac{2}{3}\ka^3\p^3T+\frac{2}{3}\ka T\p T)\ve^{(2)}\,.
$$
Thereby, we obtain the Lie algebroid ${\mathcal   A}_3(S^1)$ over $M_3(S^1)$.
Note that in ${\mathcal   A}_3(S^1)$ in contrast with the previous cases
we encounter with
the structure functions - the r.h.s of (\ref{6.8}) depends on  the projective connection $T$.

The Jacobi identuty (\ref{5.5}) in ${\mathcal   A}_3(S^1)$ takes the form
\beq{6.14}
\lb\lb\ve^{(2)}_1,\ve^{(2)}_2\rb,\ve^{(2)}_3\rb^{(1)}-
(\ve^{(2)}_1\ka\p\ve^{(2)}_2-\ve^{(2)}_2\ka\p\ve^{(2)}_1)\de_{\ve^{(2)}_3}T+{\rm c.p.}(1,2,3)
=0,
\end{equation}
\beq{6.15}
\lb\lb\ve^{(2)}_1,\ve^{(2)}_2\rb,\ve^{(1)}_3\rb^{(1)}-
(\ve^{(2)}_1\ka\p\ve^{(2)}_2-\ve^{(2)}_2\ka\p\ve^{(2)}_1)\de_{\ve^{(1)}_3}T=0.
\end{equation}
The brackets here correspond to the product of
structure functions in the left hand side
of (\ref{5.5}) and the superscript $(1)$ corresponds to the first order differential
operators.
For the  rest  brackets the Jacobi identity has the standard form.


\subsection{Global Lie algebroid over SL$(3,\mC)$-opers}

The base $M_3(\Si_{g,n}) $ of the global Lie algebroid $\clA_3(\Si_{g,n})$ are
$\SLt$-opers. They
are well defined globally on $\Si_{g,n}$. The space of its sections  $\clG(\clA_3)$  are the second order differential operators without free
terms. To define this space properly we use the formalism of Volterra operators
on  $\Si_{g,n}$ in Section 6. They are well defined on $\Si_{g,n}$. There is map
from a quotient space of the Volterra operators to  $\clG(\clA_3)$. It will be defined
in subsection 6.2.

\subsubsection{Derivation of the brackets}

To derive the Lie brackets (\ref{6.6}) - (\ref{6.8}) and the anchor action
(\ref{6.9}) -- (\ref{6.12}) globally we use the
matrix description of $\SLt$-opers (\ref{6.1}).

Consider the set $G_3(\Si_{g,n})$ of automorphisms of the bundle $E$ over $\Si_{g,n}$
\beq{6.15a}
A\to f^{-1}\kappa\p f-f^{-1}Af\,,
\end{equation}
that preserve the $\SLt$-oper structure
\beq{6.16}
f^{-1}\kappa\p f-
f^{-1}\thmat{0}{1}{0}{0}{0}{1}{W}{T}{0}f=
\thmat{0}{1}{0}{0}{0}{1}{W'}{T'}{0}\,.
\end{equation}
It is clear that $G_3(\Si_{g,n})$ is the Lie groupoid over  $M_3 (\Si_{g,n})=\{(W,T)\}$
 with
$l(f)=(W,T)$, $~r(f)=(W',T')$, $~f\to \lan\lan W,T|f|W',T'\ran\ran$.
The left identity map is
$$
P\exp(\int^z_{z_0} A(W,T))\cdot C\cdot P\exp(\int^z_{z_0} A(W,T))\,,
$$
where $C$ is an arbitrary matrix from $\SLt$ and $A(W,T))$ has the oper structure
(\ref{6.1}).
The right identity map has the same
form with $(W,T)$ replaced by $(W',T')$.

The infinitesimal  version of (\ref{6.16}) takes the form
\beq{6.17}
\kappa\p X-\left[\thmat{0}{1}{0}{0}{0}{1}{W}{T}{0},X\right]=
\thmat{0}{0}{0}{0}{0}{0}{\de W}{\de T}{0}\,.
\end{equation}
It is a linear differential system for the matrix elements of the
traceless matrix $X$. The matrix elements $x_{j,k}\in\Om^{(j-k,0)}(\Si_{g,n})$
 depend on two arbitrary fields $x_{23}=\ve^{(1)},~x_{13}=\ve^{(2)}$. The
solution takes the form
\beq{6.17b}
X=\thmat{x_{11}}{x_{12}}{\ve^{(2)}}{x_{21}}{x_{22}}
{\ve^{(1)}}{x_{31}}{x_{32}}{x_{33}}\,,
\end{equation}
$$
x_{11}=\frac{2}{3}(\ka^2\p^2-T)\ve^{(2)}-\ka\p\ve^{(1)}\,,~
x_{12}=\ve^{(1)}-\ka\p\ve^{(2)}\,,
$$
$$
x_{21}=\frac{2}{3}\ka\p(\ka^2\p^2-T)\ve^{(2)}-\ka^2\p^2\ve^{(1)}+W\ve^{(2)}\,,~
x_{22}=-\frac{1}{3}(\ka^2\p^2-T)\ve^{(2)}\,,
$$
$$
x_{31}=\frac{2}{3}\ka^2\p^2(\p^2-T)\ve^{(2)}-\ka^3\p^3\ve^{(1)}+\ka\p(W\ve^{(2)})+W\ve^{(1)}\,,
$$
$$
x_{32}=\frac{1}{3}\ka\p(\ka^2\p^2-T)\ve^{(2)}-\ka^2\p^2\ve^{(1)}+W\ve^{(2)}+T\ve^{(1)}\,,
$$
$$
x_{33}=-\frac{1}{3}(\ka^2\p^2-T)\ve^{(2)}+\ka\p\ve^{(1)}\,.
$$
The matrix elements of the commutator $[X_1,X_2]_{13}$, $[X_1,X_2]_{23}$
give rise to the brackets (\ref{6.6}), (\ref{6.7}), (\ref{6.8}). Simultaneously,
from (\ref{6.17}) one obtain the anchor action (\ref{6.9})--(\ref{6.12}).


\subsubsection{Contribution of the marked points}

Assume that the coefficients of opers $M_3=M_3(\Si_{g,n})$ have holomorphic poles at the
marked points
\beq{6.3}
T|_{z\rar x_a}\sim\frac{T^a_{-2}}{(z-x_a)^2}+\frac{T^a_{-1}}{(z-x_a)}+\ldots\,,
\end{equation}
\beq{6.4}
W|_{z\rar x_a}\sim
\frac{W^a_{-3}}{(z-x_a)^3}+\frac{W^a_{-2}}{(z-x_a)^2}+\frac{W^a_{-1}}{(z-x_a)}
+\ldots\,.
\end{equation}

The space of sections ${\mathcal   G}_3\sim\G({\mathcal   A}_3)=\{\ve^{(1)},~\ve^{(2)}\}$
of the global algebroid ${\mathcal   A}_3$  was defined above.
In addition, assume  the coefficients of
 the first and
second order differential operators
vanish holomorphically at the marked points
\beq{6.5}
\ve^{(1)}\sim r^{(1)}_a(z-x_a)+o(z-x_a)\,,~~
\ve^{(2)}\sim r^{(2)}_a(z-x_a)^2+o(z-x_a)^2\,,~~r^{(j)}\neq 0\,.
\end{equation}
Note that these asymptotics are consistent with the Lie brackets,
the anchor action and with the asymptotics (\ref{6.3}) and (\ref{6.4}).


\subsubsection{Cohomology of the Lie algebroid}

It follows from  (\ref{6.9}) -- (\ref{6.12}),  (\ref{6.3}), (\ref{6.4}),
and (\ref{6.5})
that $\de_{\ve_j}T_{-2}^a=0$, $\de_{\ve_j}W_{-3}^a=0$, $(j=1,2,\, a=1,\ldots n)$. Therefore,
\beq{co0}
T_{-2}^a\,,~W_{-3}^a\in H^0({\mathcal   A}_3,M_3)
\end{equation}
Define the cocycles
\beq{6.17a}
c^{(1)}=\int_{\Si_{g,n}}\ve^{(1)}\bp T\,,
~~c^{(2)}=\int_{\Si_{g,n}}\ve^{(2)}\bp W
\end{equation}
from $H^1({\mathcal   A}_3,M_3)$.
The contribution of the marked points to the cocycles
is equal to
$$
c^{(1)}\rar \sum_{a=1}^n r^{(1)}_aT_{-2}^a\,,~~
c^{(2)}\rar \sum_{a=1}^n r^{(2)}_aW_{-3}^a\,.
$$

The cocycles lead to the shift of the anchor action
$$
\hat{\de}_{\ve^{(j)}}f(W,T)=\lan \de_{\ve^{(j)}}W|\frac{\de f}{\de W}\ran   +
\lan \de_{\ve^{(j)}}T|\frac{\de f}{\de T}\ran   +c^{(j)}\,.
$$

There exists $2g$ central extensions $c_\al$ of  ${\mathcal   G}_3$,
provided by the nontrivial cocycles from  $H^2({\mathcal   A}_3,M_3)$. They are the
 contour integrals $\ga_\al$
\beq{6.16a}
c_\al(\ve^{(j)}_1,\ve^{(k)}_2)=\oint_{\ga_\al}\la(\ve^{(j)}_1,\ve^{(k)}_2),
~~(j,k=1,2)\,,
\end{equation}
where $\ga_\al$ are the fundamental cycles of $\Si_{g,n}$ and
$\la(\ve^{(j)}_1,\ve^{(k)}_2)$ are defined by (\ref{coc1}).
These cocycles allow us to construct the extended brackets:
$$
\lb(\ve_1^{(j)},0),(\ve_2^{(m)},0)\rb_{c.e.}=
(\lb\ve_1^{(j)},\ve_2^{(m)}\rb,\sum_\al c_\al(\ve_1^{(j)},\ve_2^{(m)}))\,,
~(j,m=1,2)\,.
$$


\subsection{Global Hamiltonian algebroid}

The affine space   ${\mathcal   R}_3$ over $T^*M_3$  is the classical phase space for the $W_3$-gravity on $\Si_g\times\mR$ \cite{P,BFK,GLM}.
Its sections are the Beltrami differentials $\mu\in\Om^{(-1,1)}(\Si_{g,n})$ and
the differentials $\rho\in\Om^{(-2,1)}(\Si_{g,n})$. They
are smooth and near the marked points behave as  (\ref{mu}) and
\beq{6.19}
\rho|_{z\rar x_a}\sim
(t^{(2)}_{a,0}+t^{(2)}_{a,1}(z-x^0_a))\bp\chi_a(z,\bz)\,.
\end{equation}

According with the general theory the anchor (\ref{6.9})--(\ref{6.12})
can be lifted from $M_3$ to ${\mathcal   R}_3$.
This lift is nontrivial owing to the cocycle (\ref{6.17a}).
It follows from (\ref{5.16a}) that the anchor action on $\mu$ and
$\rho$ takes the form
\beq{6.21}
\de_{\ve^{(1)}}\mu=-\bar{\p}\ve^{(1)}-\mu\ka\p\ve^{(1)}+\ka\p\mu\ve^{(1)}-
\rho\ka^2\p^2\ve^{(1)}\,,
\end{equation}
\beq{6.22}
\de_{\ve^{(1)}}\rho=-2\rho\ka\p\ve^{(1)}+\ka\p\rho\ve^{(1)}\,,
\end{equation}
\beq{6.23}
\de_{\ve^{(2)}}\mu=\ka^2\p^2\mu\ve^{(2)}-\frac{2}{3}\left[(\ka\p(\ka^2\p^2-T)\rho)\ve^{(2)}
-(\ka\p(\ka^2\p^2-T)\ve^{(2)})\rho\right]\,,
\end{equation}
\beq{6.24}
\de_{\ve^{(2)}}\rho=-\bar{\p}\ve^{(2)}+(\rho\ka^2\p^2\ve^{(2)}-\ka^2\p^2\rho\ve^{(2)})
+2\ka\p\mu\ve^{(2)}-\mu\ka\p\ve^{(2)}\,.
\end{equation}

The transformations (\ref{6.9}) -- (\ref{6.12}) and (\ref{6.21}) -- (\ref{6.24}) are canonical with respect to the symplectic form
$$
\om=\int_{\Si_{g,n}}D T\wedge D\mu+D W\wedge D\rho\,.
$$
They are defined by the Hamiltonians
\beq{ham5}
h^{(1)}=\int_{\Si_{g,n}}(\mu\de_{\ve^{(1)}}T+\rho\de_{\ve^{(1)}}W)+c^{(1)}\,,
\end{equation}
\beq{ham6}
h^{(2)}=\int_{\Si_{g,n}}(\mu\de_{\ve^{(2)}}T+\rho\de_{\ve^{(2)}}W)+c^{(2)}\,.
\end{equation}
After the integration by parts they take the form
$$
h^{(1)}=\int_{\Si_{g,n}}\ve^{(1)}F^{(1)}\,,~~
h^{(2)}=\int_{\Si_{g,n}}\ve^{(2)}F^{(2)}\,,
$$
where $F^{(1)}\in\Om^{(2,1)}(\Si_{g,n})$, $F^2\in\Om^{(3,1)}(\Si_{g,n})$
\beq{F1}
F^{(1)}=-\bp T-\ka^4\p^4\rho+T\ka^2\p^2\rho-\ka(3W-2\ka\p T)\p\rho-
\end{equation}
$$
-(2\ka\p W-\ka^2\p^2 T)\rho+2\ka^3\p^3\mu-2\ka\p T\mu-\ka\p T\mu\,,
$$
\beq{F2}
F^{(2)}=-\bp W-\frac{2}{5}\ka^5\p^5\rho+\frac{4}{3}T\ka^3\p^3\rho+2\ka\p T\ka^2\p^2\rho+
\ka(-\frac{2}{3}T^2+2\ka^2\p^2 T-2\p W)\p\rho
\end{equation}
$$
-(\ka^2\p^2W-\frac{2}{3}\ka^3\p^3T+\frac{2}{3}\ka T\p T)\rho+\ka^4\p^4\mu
-3W\ka\p\mu-\ka\p W\mu-\ka^2T\p^2\mu\,.
$$


\subsection{The moduli space ${\mathcal   W}_3$ of the $W_3$ gravity}

Let $G_3$ be the  groupoid corresponding to the
algebroid ${\mathcal G}_3$.
\begin{defi}
The moduli space  ${\mathcal   W}_3$ of the $W_3$-gravity is the symplectic quotient
$$
{\mathcal   W}_3={\mathcal   R}_3// G_3=\{F^1=0,F^2=0\}/G_3.
$$
\end{defi}
It has dimension
$\dim{\mathcal   W}_3=16(g-1)+6n$.
The term $6n$ comes from the coefficients
$T_{-1}^a,W^a_{-1},W^a_{-2}$, and the dual to them
$t^{(1)}_{a,0},t^{(2)}_{a,0},t^{(2)}_{a,1}$,  $(a=1,\dots,n)$ in (\ref{mu}) and
(\ref{6.19}).

The moment equations $F^{(1)}=0,~F^{(2)}=0$ are the consistency conditions
for the linear system
\beq{6.26a}
\left\{
\begin{array}{l}
(\ka^3\p^3-T\ka\p-W)\psi(z,\bz)=0\,,\\
\left(\bp +(\mu-\ka\p\rho)\p +\ka^2\rho\p^2
+\frac{2}{3}(\ka^2\p^2-T)\rho-\ka\p\mu
\right)\psi(z,\bz)=0\,,
\end{array}
\right.
\end{equation}
where $\psi(z,\bz)\in\Om^{(-1,0)}(\Si_{g,n})$.
It is an analog of the vector representation.
 We will prove this statement in next subsection.
 This system defines $W_3$-\emph{projective structure} on $\Si_{g,n}$.
In the first equation the Schrodinger operator is replaced by the
third order differential operator depending on two fields $T$ and $W$.
The second equation  represents the deformation of the
operator $\bp$ (or more general $\bp+\mu\p$ as in (\ref{8.16}))
by the second order differential
operator. The left hand side is the explicit form of the deformed operator
when it acts on the space $\Om^{-1,0}(\Si_{g,n})$. This deformation cannot
be supported by the structure of a Lie algebra  and one leaves with
the algebroid symmetries.

 Now we construct the BRST complex.  We introduce the ghosts
fields $\eta^{(1)},\eta^{(2)}$ and their momenta ${\mathcal   P}^{(1)},{\mathcal
P}^{(2)}$.
It follows from Theorem 2.1 that for
$$
\Om=\sum_{j=1,2}h^{(j)}(\eta^{(j)})+
\oh\sum_{j,k,l=1,2}\int_{\Si_{g,n}}(\lb\eta^{(j)},\eta^{(k)}\rb{\mathcal   P}^{(l)})
$$
the operator  $QF=\{F,\Om\}$ is nilpotent and define the BRST cohomology
in the complex
$$
\bigwedge{}^\bullet({\mathcal   G}_3\oplus{\mathcal   G}_3^*)\otimes C^\infty({\mathcal
R}_3)\,.
$$

\subsection{Chern-Simons derivation}

We follow here the derivation of $W$-gravity proposed in Ref.\,\cite{BFK}.
We add to this construction a contribution of the Wilson lines
 due to the presence of the marked points on $\Si_{g,n}$.

Consider the Chern-Simons functional on $\Si_{g,n}\oplus {\mR}^+$
$$
S=\int_{\Si_{g,n}\oplus \mR^+}\tr({\bf A}d{\bf A}+\frac{2}{3}{\bf A}^3)+
\sum_{a=1}^n\int_{\mR^+}\tr(p^0_a\p_tg_ag_a^{-1})\,,~~
({\bf A}=(A,\bA,A_t))\,,
$$
where the last sum is the geometric action coming from the Kirillov-Kostant forms on the
coadjoint
orbits ${\mathcal O}_a$ (\ref{orb}).
 Introduce $n$ Wilson lines $W_a(A_t)$
along the time directions and located at the marked points
$$
W_a(A_t)=P\exp \tr(p_a\int A_t),~a=1,\ldots,n\,.
$$
In the hamiltonian picture the phase space, corresponding to the Chern-Simons functional is
\beq{psg}
{\mathcal R}_{SL_3}=
\{A\,,\bA\,,{\mathcal O}_1,\ldots,{\mathcal O}_n\}\,,
\end{equation}
endowed with the symplectic form
(\ref{8.4a}).  The field $A_t$ is the Lagrange multiplier for the first class constraints
(\ref{8.5a}).

The phase space of $W_3$-gravity ${\mathcal   R}_3$ can be derived from ${\mathcal
R}_{SL_3}$.
 The flatness condition (\ref{8.5a}) generates
the gauge transformations
\beq{5.31}
A\to f^{-1}\ka\p f-f^{-1}Af\,,~~\bA\to f^{-1}\bp f-f^{-1}\bA f\,,~~
g_a\to g_af_a\,.
\end{equation}
The symplectic quotient with respect to the gauge group $G_{\SLt}$
is the moduli
space ${\mathcal   M}^{flat}_3$ of the flat $\SLt$ bundles over $\Si_{g,n}$.

Let $P$ be the maximal parabolic subgroup of $\SLt$ of the form
$$
P=\thmat{*}{*}{0}{*}{*}{0}{*}{*}{*}\,,
$$
and $G_P$ be the corresponding gauge group.
We partly fix first the gauge with respect to $G_P$.
A generic connection $\nabla$ can be transformed by $f\in G_P$ to
  the form (\ref{6.1}).
Taking into account (\ref{8.5a}) we assume that $A$ has simple poles at the marked points.
To come to $M_3$
one should respect the behavior of the matrix elements at the marked points
(\ref{6.3}), (\ref{6.4}). For this purpose we use an additional singular gauge transform
by the diagonal matrix
$$
h=\prod_{a=1}^n\chi_a(z,\bz)\di((z-x_a)^{-1},1,(z-x_a))\,.
$$
The resulting gauge group we denote $G_{(P\cdot h)}$, where $\chi_a$ is defined by
(\ref{cf}).

The form of $\bA$ can be read off from (\ref{8.5a})
\beq{5.32}
\bA=\thmat{a_{11}}
{a_{12}}
{-\rho}
{a_{21}}
{a_{22}}
{-\mu}
{a_{31}}
{a_{32}}
{a_{33}}\,,
\end{equation}
$$
a_{11}=-\frac{2}{3}(\ka^2\p^2-T)\rho+\ka\p\mu\,,~~
a_{12}=-\mu+\ka\p\rho\,,
$$
$$
a_{21}=-\frac{2}{3}\ka\p(\ka^2\p^2-T)\rho+\ka^2\p^2\mu-W\rho\,,~~
a_{22}=\frac{1}{3}(\ka^2\p^2-T)\rho\,,
$$
$$
a_{31}=-\frac{2}{3}\ka^2\p^2(\ka^2\p^2-T)\rho+\ka^3\p^3\mu-\ka\p(W\rho)-W\mu\,,
$$
$$
a_{32}=-\frac{1}{3}\ka\p(\ka^2\p^2-T)\rho+\ka^2\p^2\mu-W\rho-T\mu\,,~~
a_{33}=\frac{1}{3}(\ka^2\p^2-T)\rho-\ka\p\mu\,.
$$
The condition (\ref{8.5a})  for the special choice $A$ (\ref{6.1}) and
$\bA$ (\ref{5.32}) gives rise to the relations
$F(A,\bA)|_{(3,1)}=F^{(2)}$ (\ref{F1}), $F(A,\bA)|_{(2,1)}=F^{(1)}$ (\ref{F2}),
while the other matrix
elements of $F(A,\bA)$ vanish identically.
At the same time, the matrix linear
system (\ref{8.5c}) coincides with  (\ref{6.26a}).
In this way, we come to the matrix description of the moduli space
${\mathcal   W}_3$.

The cocycles $c_\al(\ve^{(j)}_1,\ve^{(k)}_2)$ (\ref{6.16a}) can be derived from
the two-cocycle (\ref{8.3}) of ${\mathcal   A}_{SL_3}$. Substituting in (\ref{8.3})
the matrix realization of $\G({\mathcal   A}_3)$ (\ref{6.17b}), one comes to (\ref{6.16a}).

The action of groupoid $G_3$ on $A,\bA$ plays the role of the rest gauge
transformations that complete the $G_{(P\cdot h)}$ action to the $G_{SL_3}$ action.
The algebroid symmetry with non-trivial structure functions arises in this theory as a result of the partial
gauge fixing by $G_{(P\cdot h)}$. Thus we come to the following diagram

\bigskip
$$
\begin{array}{rcccl}
           &\fbox{${\mathcal   R}_{SL_3}$}&                       &                    &\\
           &       |                      &\searrow{G_{(P\cdot h)}}&           &\\
G_{\SLt} &        |                     &                    &\fbox{${\mathcal   R}_3$}&\\
           &\downarrow    &                   &\downarrow &G_3\\
           &\fbox{${\mathcal   M}^{flat}_{SL_3}$}&       &\fbox{${\mathcal   W}_3$}&\\
\end{array}
$$
\bigskip

The tangent space to ${\mathcal   M}^{flat}_{SL_3}$ at the point $A=0$, $\,\bA=0$,
$\,p_a=0$, $\,g_a=id$
coincides with the tangent space to ${\mathcal   W}_3$ at the point
$W=0$, $\,T=0$, $\,\mu=0$, $\,\rho=0$. Their dimension is $16(g-1)+6n$. But their global
structure is different and the diagram cannot be closed by the horizontal
isomorphisms. The interrelations between ${\mathcal   M}^{flat}_{SL_N}$ and ${\mathcal
W}_N$
were analyzed in \cite{H2,Go}.


\section{AGD algebroids and generalized  projective structures}
\setcounter{equation}{0}

In this section we define  generalized projective structures on $\Si_{g,n}$
 related the  $\GLN$ and  $\SLN$-opers. It is a phase space of $W_N$-gravity.
In particular, we define
deformations of complex structures
(the $W_N$-deformations)
by the Volterra operators and by the opers. To construct a Lie algebroid over
 the space of $\GLN$-opers we use  the pairing (\ref{pa1}) corresponding to the case (\ref{tab2}) with $m=N$. As a result the space of sections of the Lie algebroid
is the space the Volterra operators instead of the space the differential operators, considered in previous Section.
We start with the description of the local AGD algebroid following Ref.\,\cite{GD}
 and then give its global version. The passage from the Lie
algebroid to the Hamiltonian algebroid allows us to describe
the generalized projective structures and their moduli.

\bigskip

\subsection{Local AGD algebroid }

Consider a set $B=\Psi DO(D)$ of pseudo-differential operators on a  disk $D\subset \Si_{g,n}$ and their restriction on the boundary $S^1\sim\p D$.
It is a ring of formal Laurent series
$$
B=B((\p^{-1}))=
\{B_{r,N}(S^1)\,,~~r,N\in\mZ\}=\{X(t,\p)\}\,,~t\in S^1\,,~\p=\p_t
$$
\beq{9.1}
X(t,\p)=\sum_{k=-\infty}^{r-1}a_k(t)\p^k,~~
( a_k(z)\in \Om^{-N-k+1}(S^1) \,.
\end{equation}
 The multiplication on $B$ is defined as the non-commutative
multiplication of their symbols
\beq{pr}
X(t,\la)\circ Y(t,\la)=
\sum_{k\geq 0}\f1{k!}\frac{\p^k}{\p\la^k} X(t,\la)\frac{\p^k}{\p t^k}Y(t,\la)\,.
\end{equation}
In what follows we omit the multiplication symbol $\circ$.

Note that $B_{r,N}(S^1)\in\Psi DO(S^1)$ can be considered as the formal map
of the sheaves
\beq{B}
B_{r,N}(S^1)~:~\Om^{\frac{N-1}{2}}(S^1)\to\Om^{-\frac{N-1}{2}}(S^1)\,.
\end{equation}

Let  $M^G_N(S^1)$ be a space of differential operators $L_N=\ka^N\p^N+W_1\ka^{N-1} \p^{N-1}+\ldots+W_N$  on $S^1\,$
 $\,(\p=\p_t)$ with smooth coefficients, corresponding to the $\GLN$-oper on $D$.
 For brevity we call them the $\GLN$-opers.
Then we have
$$
XL_N~:~\Om^{-\frac{N-1}{2}}(S^1)\to\Om^{-\frac{N-1}{2}}(S^1)\,,~~~
L_NX~:~\Om^{\frac{N+1}{2}}(S^1)\to\Om^{\frac{N+1}{2}}(S^1)\,,
$$
where the product is defined by (\ref{pr}).

Define a pairing between $M^G_N(S^1)$  and $B_{r,N}(S^1)$.
 For $L_NX=\sum_k c_k\p^k$  let $Res\,L_NX=c_{-1}$. It is a one-form on
 $S^1$ and one can define a pairing
\beq{pa1}
\lan L_NX\ran=\f1{2\pi}\oint_{S^1} Res (L_NX)dt\,.
\end{equation}
  Because $Res\,[L_N,X]$ is a derivative, $\lan L_NX\ran=\lan XL_N\ran$. The pairing corresponds to the case (\ref{tab2}).

\bigskip

\subsubsection{Local AGD algebroid over $\GLN$-opers}

The local AGD algebroid over $\GLN$-opers was constructed implicitly in \cite{GD}.
The AGD brackets on the space $M^G_N(S^1)$
are defined as follows.
The space of sections of the cotangent bundle $T^*M^G_N(S^1)$ can be identified
with  the quotient space of the Volterra operators
\beq{9.1a}
\G(T^*M^G_N(S^1))= B_{0,N}(S^1)/B_{-N-2,N}(S^1)\,.
\end{equation}
As before consider $\ka$ as an independent variable of an oper $L_N$ and
extend the pairing (\ref{pa1}) by  the additional term $\ka\cdot c$, where $c$  corresponds to the central extension. In this case instead of (ref{9.1a})
we come to the extended space of sections
$$
\hat{\G}(T^*M^G_N(S^1))= B_{0,N}(S^1)/B_{-N-2,N}(S^1)\oplus \mC\,.
$$

For $L_N $ and $X\in \G(T^*M^G_N(S^1))$ define the functional
$l_X=\lan L_NX\ran$.
In particular, for
$l_X=W_1(z)\,$,  $\,X=\de(z/t)\p_t^{-N}$.
The AGD brackets have the form
\beq{9.2}
\{l_X,l_Y\}=\lan L_NX(L_NY)_+\ran -\lan XL_N(YL_N)_+\ran\,,
\end{equation}
where $X_+=\sum_{k=0}^Na_k\p^k $ is the differential part of $X$.
Equivalently, (\ref{9.2}) can be rewritten as
\beq{9.2c}
\{l_X,l_Y\}= \lan XL_N(YL_N)_-\ran-\lan L_NX(L_NY)_-\ran \,,
\end{equation}
 where $A_-$ is the integral part of $A\in\Psi DO\,$, $\,(A_+=A-A_-)$.

Using the general prescription (\ref{2.3}) we find from (\ref{9.2c})
the Lie brackets in the space  of sections $\hat{\G}(T^*M^G_N(S^1))$
\beq{9.3}
\lb (X,0),(Y,0)\rb=\left(((YL_N)_-X-X(L_NY)_-+X(L_NY)_+-(Y(L_N)_+X),c\right)\,,
\end{equation}
where
\beq{centr}
c=\frac{\p}{\p\ka}\{l_X,l_Y\}\,.
\end{equation}
Due to the Jacobi identity for the brackets (\ref{9.2}) this term leads to
the central extension of the Lie algebra $\G(T^*M^G_N(S^1))$.

The anchor map assumes the form (see (\ref{7.40}))
\beq{9.4}
\de_YL_N=(L_NY)_+L_N-L_N(YL_N)_+\,.
\end{equation}
\begin{defi}
The AGD algebroid $\clA^G_N(S^1)$ over $M^G_N(S^1)$ is a bundle $ T^*M^G_N(S^1)$
with the brackets (\ref{9.3}) and the anchor (\ref{9.4}).
\end{defi}

It follows from (\ref{9.2}) and (\ref{9.4}) that the Poisson brackets can be rewritten as
\beq{9.4a}
\{l_X,l_Y\}=\lan X\de_YL_N\ran\,,
\end{equation}
or in the form of the "Poisson-Lie brackets"
\beq{9.4b}
\{l_X,l_Y\}=\oh\lan \lb X,Y\rb L_N\ran= l_{\oh\lb X,Y\rb}\,.
\end{equation}
The coefficient $1/2$ arises from the quadratic form of the Poisson bivector.
These two representations implies that the anchor plays the role of the coadjoint
action.
\begin{rem}
It is assumed in  (\ref{9.3}) and (\ref{9.4}) that the sections $X$ and $Y$ are independent on 
a point in the base $M^G_N(S^1)$. To pass to a generic section one should use the defining
properties of Lie algebroids. 
In fact, the Lie brackets (\ref{9.3}) and the anchor (\ref{9.4}) were appeared already in
\cite{GD}. 
\end{rem}

\subsubsection{Local AGD Lie algebroid over $\SLN$-opers}

 Now we consider the space $M_N(S^1)$ of $\SLN$-opers   and construct the corresponding ${\mathcal   A}_N(S^1)$ algebroid. Remind that
 an $\SLN$-oper $L^S_N$ on $D$ is defined by the condition $W_1=0\,$,
$\,(L^S_N=\p^N+W_2\p^{N-2}+\ldots)$.

The description of the AGD $\SLN$-algebroids is based on the following
statement
\begin{predl}
\begin{itemize}

  \item The anchor action (\ref{9.4}) preserves the coefficient $W_1$ of
  a $\GLN$-oper  iff
\beq{sop}
 Res\, [Y,L_N]=0\,.
\end{equation}
In particular it  preserves an $\SLN$-oper $L_N^S$.
  \item Sections satisfying (\ref{sop}) generate a Lie algebra $\clG$.
\end{itemize}
\end{predl}

This statement leads to the following definition.
\begin{defi}
 The AGD Lie algebroid ${\mathcal   A}_N(S^1)$ over the
$\SLN$-opers on $S^1$ is defined by the brackets  (\ref{9.3}) and the anchor map (\ref{9.4}) with condition (\ref{sop}).
\end{defi}

\emph{Proof of Proposition}.
If $\de_YL^S_N$ does not change $W_1$ then $\lan a\p_t^{-N}\de_YL_N\ran=0$
for any continues functional $a$ on the space $C^\infty(S^1)$.
Due to (\ref{9.4}) it implies
$$
\lan a\p_t^{-N}(L_NY)_+L_N-a\p_t^{-N}L_N(YL_N)_+\ran=0\,.
$$
 We rewrite the l.h.s. as
$$
\lan a\p_t^{-N}L_N(YL^S_N)_--a\p_t^{-N}(L_NY)_-L_N\ran
$$
$$
=\lan (a\p_t^{-N}L_N)_+YL_N-(L_Na\p_t^{-N})_+L_NY=\lan a[Y,L_N]\ran\,.
$$
Vanishing of this expression  for any $a$ is equivalent to (\ref{sop}).

Let $\de_X$ and $\de_Y$ preserve the structure of the $\SLN$-oper.
Since $[\de_X,\de_Y]=\de_{\lb X,Y\rb}$ the sections, satisfying  (\ref{sop})
generate a Lie algebra $\clG$.
Moreover,  $\SLN$-opers generate a Poisson subalgebra. In fact, we have from
(\ref{9.4a}) and (\ref{9.4b})
$$
\de_{\lb X,Y\rb}W_1=\{l_{\lb X,Y\rb},W_1\}=2\{\{l_X,l_Y\},W_1\}=0\,.
$$
$\Box$

Now prove that $\clG$ is isomorphic to $\hat{\G}(T^*M_N(S^1))\subset\hat{\G}(T^*M^G_N(S^1))$.
Consider the cotangent bundle $ T^*M_N(S^1)$ to the space of $\SLN$-opers
$M_N(S^1)$. The space of its sections is quotient space of
$T^*M_N(S^1)= T^*M^G_N(S^1)/\{a(t)\p_t^{-N}\}$.
It is possible to choose a section $X$ of the $\SLN$ Lie algebroid
 ${\mathcal   A}_N(S^1)$
such that $X\in\G(T^*M_N(S^1)$.
\begin{lem}
For any $X\in\G( T^*M^G_N(S^1))$ one can find $a(t)\p_t^{-N}$ such that
$Y=X+a(t)\p_t^{-N}$ obeys (\ref{sop}).
\end{lem}
\emph{Proof.}
It easy to find that for  $\SLN$-oper  $L_N^S\,$ $\,Res\,[a\p_t^{-N},L_N^S]=-N\p_ta$.
 For any $X\in\G(T^*M^G_N(S^1)\,$ $\,Res\,[X,L^S_N]=\p_tF(t)$.
Then $Res\, [(X+a\p_t^{-N},L^S_N]=\p_t(F(t)-Na(t))$. Choosing $a(t)=\f1{N}F(t)$
we obtain a section of the $\SLN$ algebroid. $\Box$


\subsection{Global AGD Lie algebroid}

As it was mentioned above the opers are well defined globally on the curves.
We assume that in  neighborhoods  of the marked points
the coefficients $W_j$ behave as
\beq{6.4a}
W_j|_{z\rar x_a}\sim
W^a_{-j}(j)(z-x_a)^{-j}+W^a_{-j}(j-1){(z-x_a)^{-j+1}}+
\ldots\,.
\end{equation}

The base of the global AGD $\GLN$ algebroid ${\mathcal   A}^G_N(\Si_{g,n})$ (${\mathcal   A}_N(\Si_{g,n})$) is the space  $M^G_N(\Si_{g,n})$ of global $\GLN$-opers on $\Si_{g,n}$. Similarly,
the base  of the global AGD $\SLN$ algebroid ${\mathcal   A}_N(\Si_{g,n})$  is the space  $M_N(\Si_{g,n})$ of global $\SLN$-opers.
They have the prescribed behavior
near the marked points. The spaces of their sections $\clG_N^G\,$ ($\clG_N$)
 are the quotient spaces of the Volterra operators
$$
\clG_N^G=\G({\mathcal   A}^G_N(\Si_{g,n}))= B_{0,N}(\Si_{g,n})/B_{-N-2,N}(\Si_{g,n})\,,
$$
where
$$
B_{r,N}(S^1)~:~\Om^{\frac{N+1}{2}}(\Si_{g,n})\to\Om^{-\frac{N-1}{2}}(\Si_{g,n})\,.
$$
Near a marked point $x_a$ with a local coordinate $z$ a section $X\in {\mathcal   G}_{N}$
has the expansion
\beq{X}
X=\sum_{j=1}^{N+1}\ep^{(j)}\p_z^{-j}\,,
\end{equation}
where
$\ep^{(j)} \sim r_a ^{(j)}(z-x_a)^j+o(z-x_a)^j\,,~~ r_a ^{(j)}\neq 0\,.$

Similarly,
the base  of the global AGD $\SLN$ algebroid ${\mathcal   A}_N(\Si_{g,n})$  is the space  $M_N(\Si_{g,n})$ of global $\SLN$-opers. The space of sections $\clG_N^G$ of ${\mathcal A}_N(\Si_{g,n})$  satisfy (\ref{sop}).

Consider in detail the case $N=3$. Locally the sections of $\clG_3^G$ can be represented
by the operators
$$
X=\ep^{(1)}\p_z^{-1}+\ep^{(2)}\p_z^{-2}+\ep^{(3)}\p_z^{-3}+\p_z^{-4}\,.
$$
The space of sections described in Sections 5.2 and 5.3 (\ref{m3}) are the second order
differential operators $\{ (\ve^{(1)}\frac{d}{dz}+ \ve^{(2)}\frac{d^2}{dz^2})\}$.
We express $\ve^{(1)}$ and $\ve^{(2)}$ in terms of $\ep^{(j)}$.
Note that for SL$(3,\mC)$-opers (\ref{sop}) takes the form
$$
\p_z((\p_z^2-T)\ep^{(1)}+3\ep^{(2)}+3\ep^{(3)})=0\,.
$$
 Then the anchor (\ref{9.4}) and the brackets
(\ref{9.3}) coincide with the anchor (\ref{6.9}) -- (\ref{6.12})
and the brackets (\ref{6.6}) -- (\ref{6.8}), if one puts $\ep^{(1)}= \ve^{(2)}$
and $\ep^{(2)}=\ve^{(1)}$.

The one-cocycle representing $H^1({\mathcal   A}_N)$ comes from
 the integration over $\Si_{g,n}$
$$
c(L_N,X)=\int_{\Si_{g,n}}Res(X{\bp}L_N)\,.
$$
The contribution of the marked points is
$$
c(L_N,X)|_{z=x_a}= \sum_{j=1}^Nr_a ^{(j)}W_{-j}^a(j)\,.
$$

Let $\ga_\al$ be a set of $2g$ fundamental cycles of $\Si_{g,n}$.
One can define local AGD brackets $\{l_X,l_Y\}_\al$ (\ref{9.2}) by the pairing (\ref{pa1}) using $\ga_\al$. In this way we obtain $2g$ generators of
$H^2({\mathcal   A}_N)$ of type (\ref{centr}).

\subsection{Global Hamiltonian AGD algebroid}

Let
$B_{(r,N,1)}( \Si_{g,n})= B_{r,N}(\Si_{g,n})\otimes \bar{K}(\Si_{g,n})\,,$,
where $\bar{K}$ is the anti-canonical class,
$$
B_{(r,N,1)}( \Si_{g,n})\,:
\,\G(\Om^{(\frac{N+1}{2},0)} ( \Si_{g,n}))
\to\G(\Om^{(\frac{-N-1}{2},1)} ( \Si_{g,n}))\,.
$$
Consider the set of the $\SLN$-opers $M_N=M_N( \Si_{g,n})$.
The affine space ${\mathcal   R}_N=Aff\,T^*M_N( \Si_{g,n})$
over $T^*M_N( \Si_{g,n})$
is the set of fields \\
$\xi\in B_{(0,N,1)}( \Si_{g,n})/ B_{(-N-2,N,1)}( \Si_{g,n})$.
Near the marked points $\xi$ behaves as
$$
\xi=\sum_{j=1}^{N+1}\nu_j\p^{-j},~~(\nu_{N+1}=1)\,,
$$
$$
\nu_j\sim
(t^{(j)}_{a,0}+\ldots+ t^{(j)}_{a,j-1}(z-x_a)^{j-1})\bar{\p}\chi_a(z,\bz)\,,
~~(\nu_j\in \Om^{(j-N,1)}(\Si_{g,n}))\,.
$$
The symplectic form on  ${\mathcal   R}_N$ is
$$
\om= \int_{\Si_{g,n}}Res\,(DL_N \wedge D\xi)\,.
$$
The anchor action on $M_N$ (\ref{9.4}) can be lifted from $M_N$
to ${\mathcal   R}_N$ as the canonical transformations of $\om$
\beq{9.7}
\de_Y\xi=-\bp Y+
Y(L_N\xi)_+-(\xi L_N)_+Y+(YL_N)_+\xi-\xi(L_NY)_+\,.
\end{equation}

 The transformations are generated by the Hamiltonians
$$
h_Y=\int_{\Si_{g,n}}Res\,(\xi\de_{Y}L_N)+c(L_N,Y)\,.
$$
The anchor action
$$
\de_YL_N=\{h_Y,L_N\}\,,~~\de_Y\xi=\{ h_Y,\xi\}\,.
$$
defines the global Hamiltonian AGD-algebroid ${\mathcal A}^H_N(\Si_{g,n})$.
The Hamiltonian can be represented in the form
$$
h_Y=\int_{\Si_{g,n}}Res\,(YF(L_N,\xi))\,,
$$
where
\beq{9.9}
F:=\bp L_N-(L_N\xi)_+L_N+L_N(\xi L_N)_+\,.
\end{equation}

The space $\clR_N$ is a phase space of $W_N$-gravity in the space $\Si_g\times\mR$.
The canonical transformations (\ref{9.4}), (\ref{9.7}) are the gauge transformations
of the theory.


\subsection{Generalized projective structures}

Define the operator
\beq{9.11}
\bp+A\,:\,
\Om^{(\frac{N+1}{2},0) }(\Si_{g,n})\to\Om^{(\frac{N+1}{2},1)} (\Si_{g,n})\,,~~A=-(L_N\xi)_+
\end{equation}
and the dual operator
$$
\bp+A^*\,:\,
\Om^{(-\frac{N-1}{2},0)} (\Si_{g,n})\to\Om^{(-\frac{N-1}{2},1)} (\Si_{g,n})\,,
~~A^*=(\xi L_N)_+\,.
$$
Let $\psi=(\psi^-,\psi^+)$, $\psi^-\in\Om^{(-\frac{N-1}{2},0)} (\Si_{g,n})$,
$ \psi^+\in\Om^{(\frac{N+1}{2},0)} (\Si_{g,n})$.
 Similar to Lemma 3.2 we find that the constraints $F=0$
(see (\ref{9.9})) are equivalent to the linear problem
\beq{9.10}
L_N\psi^-(z,\bz)=0\,,
\end{equation}
\beq{9.10a}
(\bp-(L_N\xi)_+)\psi^-(z,\bz)=0\,,
\end{equation}
\beq{9.12}
(\bp+(\xi L_N)_+)\psi^+(z,\bz)=0\,.
\end{equation}
Here we use the "vector representation". It means that in the matrix forms of
opers (\ref{opN}) $\psi^\pm$ are vectors. The linear system defines \emph{the generalized projective structures} on $\Si_{g,n}$.

In this way an oper $L_N$ together with the dual element $\xi$
defines $W_N$-{\em deformation of  complex structures} on $\Si_{g,n}$.
The equations (\ref{9.10a}) and (\ref{9.12}) are equivalent to the deformed
holomorphity condition for the sections\\
$\Om^{(-\frac{N-1}{2},0)} (\Si_{g,n})$
and $\Om^{(\frac{N+1}{2},0)} (\Si_{g,n})$.

Let $G_N$ be the groupoid corresponding to the AGD-algebroid ${\mathcal A}_N(\Si_{g,n})$.
\begin{defi}
The moduli space ${\mathcal W}_N$ of the $W_N$-gravity is the symplectic quotient
$$
\clR_N//G_N=\{\bp L_N-(L_N\xi)_+L_N+L_N(\xi L_N)_+=0\}/G_N\,.
$$
\end{defi}
The moduli space of the $W_N$-deformations of complex structures
on $\Si_{g,n}$ is a part of the symplectic quotient
${\mathcal   W}_N\sim{\mathcal   R}_N//G_{N}$.
The cohomology of the classical BRST operator are defined by
$$
\Om=h_\eta+
\oh\int_{\Si_{g,n}}Res\left(
(\lb\eta,\eta'\rb{\mathcal   P}
\right)\,,
$$
where $\eta$ is the ghost field corresponding to the gauge field $Y$
and ${\mathcal   P}$ is its momenta.

\small{

}
\end{document}